\documentclass[12pt]{article}
\usepackage{graphicx}
\usepackage{amsmath,amssymb,amsthm,amsfonts}
\newtheorem{thm}{Theorem}[section]
\newtheorem{cor}[thm]{Corollary}
\newtheorem{lem}[thm]{Lemma}
\newtheorem{prop}[thm]{Proposition}

\newcommand{\lam}{\lambda}

\newcommand{\R}{{\mathbb{R}}}
\newcommand{\Z}{{\mathbb{Z}}}

\newcommand{\1}{\partial}
\newcommand{\2}{\overline}
\newcommand{\3}{\varepsilon}
\newcommand{\4}{\widetilde}

\begin{document}
\title{Existence and dynamic properties of\\ 
a parabolic nonlocal MEMS equation}
\author{Kin Ming Hui\\
Institute of Mathematics, Academia Sinica,\\
Nankang, Taipei, 11529, Taiwan, R. O. C.\\
e-mail:kmhui@gate.sinica.edu.tw}
\date{Aug 16, 2010}
\smallbreak \maketitle
\begin{abstract}
Let $\Omega\subset\mathbb{R}^n$ be a $C^2$ bounded domain and $\chi>0$
be a constant. We will prove the existence of constants 
$\lambda_N\ge\lambda_N^{\ast}\ge\lambda^{\ast}(1+\chi\int_{\Omega}
\frac{dx}{1-w_{\ast}})^2$ for the nonlocal MEMS equation $-\Delta v
=\lam/(1-v)^2(1+\chi\int_{\Omega}1/(1-v)dx)^2$ in $\Omega$, $v=0$ 
on $\1\Omega$, such that a solution exists for any $0\le\lambda
<\lambda_N^{\ast}$ and no solution exists for any $\lambda
>\lambda_N$ where $\lambda^{\ast}$ is the pull-in voltage 
and $w_{\ast}$ is the limit of the minimal solution of 
$-\Delta v=\lam/(1-v)^2$ in $\Omega$ with $v=0$ on $\1\Omega$
as $\lambda\nearrow
\lambda^{\ast}$. Moreover $\lambda_N<\infty$ if 
$\Omega$ is a strictly convex smooth bounded domain. We will prove the 
local existence and uniqueness of the parabolic nonlocal 
MEMS equation $u_t=\Delta u+\lam/(1-u)^2(1+\chi\int_{\Omega}1/(1-u)\,dx)^2$ 
in $\Omega\times (0,\infty)$, $u=0$ on $\1\Omega\times (0,\infty)$, 
$u(x,0)=u_0$ in $\Omega$. We prove the existence of a unique global solution 
and the asymptotic behaviour of the global solution of the parabolic 
nonlocal MEMS equation under various boundedness conditions
on $\lambda$. We also obtain the quenching 
behaviour of the solution of the parabolic nonlocal MEMS equation when 
$\lambda$ is large.  
\end{abstract}

\vskip 0.2truein


\vskip 0.2truein
\setcounter{equation}{0}
\setcounter{section}{-1}

\section{Introduction}
\setcounter{equation}{0}
\setcounter{thm}{0}

Micro-electromechanical systems (MEMS) are widely used nowadays in many 
electronic devices including accelerometers for airbag deployment in cars, 
inkjet printer heads, and the device for the protection of hard disk, etc. 
The challenge is to build and understand the mathematical models and the 
mechanism for the various MEMS devices. Recently there is a lot of study 
on the equations arising from MEMS by P.~Esposito, N.~Ghoussoub, 
Y.~Guo, Z.~Pan and M.J.~Ward \cite{EGG1}, \cite{GG1}, \cite{GG2}, 
\cite{GPW}, N.I.~Kavallaris, 
T.~Miyasita and T.~Suzuki \cite{KMS}, F.~Lin and Y.~Yang \cite{LY},
L.~Ma and J.C.~Wei \cite{MW}, G. Flores, G.A.~Mercado, J.A.~Pelesko
and A.A.~Triolo  \cite{FMP}, \cite{P}, \cite{PT} etc. Interested readers 
can read the book, ``Modeling MEMS and NEMS'' \cite{PB}, by J.A.~Pelesko and 
D.H.~Bernstein for the mathematical modeling and various applications of 
MEMS devices.
 
In \cite{PB} J.A.~Pelesko and D.H.~Berstein model the deflection between
the two parallel plates of an electrostatic MEMS device by the 
equation
\begin{equation*}\displaystyle
\left\{\begin{aligned}
-\Delta w=&\frac{\lam}{(1-w)^2}\quad\mbox{ in }\Omega\\
w=&0\qquad\qquad\mbox{on }\1\Omega
\end{aligned}\right.\tag {$S_\lambda$}
\end{equation*}
where $\Omega\subset\R^2$ is a bounded $C^2$ domain.
Interested readers can read the papers \cite{GG1}, \cite{KMS}
and \cite{LY} for various results on the above equation.
In \cite{LY} F.H.~Lin and Y.~Yang by using 
variational argument derived the following nonlocal MEMS equation
\begin{equation*}\displaystyle
\left\{\begin{aligned}
-\Delta v=&\frac{\lam}{(1-v)^2(1+\chi\int_{\Omega}\frac{dx}{1-v})^2}
\quad\mbox{ in }\Omega\\
v=&0\quad\qquad\qquad\qquad\qquad\quad\,\,\mbox{on }\1\Omega
\end{aligned}\right.\tag {$S^N_\lambda$}
\end{equation*} 
of an electrostatic MEMS device with circuit series capacitance
that models the deflection between a membrane and an upper plate which 
is parallel to the plane containing the boundary of the membrane. 
An interesting 
property of ($S_{\lambda}$) (\cite{GG1}, \cite{LY}) is that there exists 
$\lambda^{\ast}>0$ called pull-in voltage in the literature of MEMS 
research such that ($S_{\lambda}$) has a solution
for any $0\le\lambda<\lambda^{\ast}$ and no solution exists for any $\lambda
>\lambda^{\ast}$. Physically this corresponds to the existence of a 
pull-in voltage such that the membrane and the upper plate 
in the MEMS device collapse together \cite{LY}, \cite{PB}, when $\lambda$ 
which is proportional to the square of the electric voltage between 
the membrane and the upper plate 
is greater than the pull-in voltage $\lambda^{\ast}$. 

In this paper we will study the equation ($S^N_\lambda$) and show that
($S^N_\lambda$) has similar properties. Let $\chi>0$. We will study the 
existence and non-existence of solutions of the corresponding nonlocal 
parabolic equation (cf. \cite{PB}, \cite{EGG2}), 
\begin{equation*}\left\{\begin{aligned}
\frac{\1 u}{\1 t}=&\Delta u
+\frac{\lambda}{(1-u)^2(1+\chi\int_\Omega\frac{dy}{1-u(y,t)})^2} 
\quad\quad\mbox{in }\Omega\times (0,T)\\
u=&0\qquad\qquad\qquad \qquad\qquad\qquad\quad\quad\quad\quad
\mbox{ on }\partial\Omega\times (0,T)\\
u(x,0)=&u_0 \qquad\quad\qquad\qquad\qquad\qquad\qquad\quad\quad\,\,
\mbox{ in }\Omega
\end{aligned}\right.\tag {$P_{\lambda}$}
\end{equation*}
where $\lambda\ge 0$ is a constant. The above equation also appears in the
unpublished preprint ``Pull-in voltage and steady states of nonlocal 
electrostatic MEMS'' of N.~Ghoussoub and Y.~Guo.
We will prove the local existence and 
uniqueness of solution of ($P_{\lambda}$). Under some boundedness conditions
for $\lambda$ we prove the existence of a unique global solution and the
asymptotic behaviour of the global solution of ($P_{\lambda}$). We prove 
the quenching behaviour of the solution of ($P_{\lambda}$) when 
$u_0\equiv 0$ on $\Omega$ and $\lambda$ is large. Physically this 
corresponds to the case that there is no deflection of the plates at 
the initial time $t=0$ and the applied voltage is large. We also prove
the quenching behaviour of the solution of ($P_{\lambda}$)
when $\Omega$ is a ball, $u_0$ is radially symmetric, and $\lambda$ is large.  

The plan of the paper is as follows. In section 1 we will prove the existence
of constants $\lambda_N\ge\lambda_N^{\ast}\ge\lambda^{\ast}(1+\chi\int_{\Omega}
\frac{dx}{1-w_{\ast}})^2$ such that ($S_{\lambda}^N$) 
has a solution for any $0\le\lambda<\lambda_N^{\ast}$ and ($S_{\lambda}^N$) 
has no solution for any $\lambda>\lambda_N$. We also prove the boundedness 
of $\lambda_N$ when 
$\Omega$ is a strictly convex smooth bounded domain of $\R^n$. In 
section 2 we will prove the local existence and uniqueness of solution 
of ($P_{\lambda}$). We also obtain energy estimates for
the solution of ($P_{\lambda}$). In section 3 we prove the global 
existence and asymptotic behaviour of the global solution of ($P_{\lambda}$) 
under various boundedness conditions on $\lambda$. In section 4 
we prove the quenching behaviour of the solution 
of ($P_{\lambda}$) when $\lambda$ is large.

We will assume that $\Omega\subset\R^n$ is a bounded $C^2$ domain for the 
rest of the paper. We start with some definitions. For any $\delta>0$, 
$R>0$, let $\Omega_{\delta}=\{x\in\Omega:\mbox{dist}(x,\1\Omega)<\delta\}$ 
and $B_R=\{x\in\R^n:|x|<R\}$. We say that $w$ is a 
solution of ($S_{\lambda}$) (($S_{\lambda}^N$) 
respectively) if $w\in C^2(\Omega)\cap C(\2{\Omega})$, $0\le w<1$ in 
$\Omega$, satisfies ($S_{\lambda}$) (($S_{\lambda}^N$) respectively) 
in the classical sense. 

For any constants $\chi\ge 0$, $\lambda>0$, $f\in C(\2{\Omega}\times 
(0,T))$ and
\begin{equation}
u_0\in L^1(\Omega)\text{ with }u_0\le a\text{ a.e. in }\Omega
\end{equation}
for some constant $0<a<1$ we say that $u$ is a solution 
(subsolution, supersolution respectively) of 
\begin{equation*}\left\{\begin{aligned}
\frac{\1 u}{\1 t}=&\Delta u
+\frac{\lambda f}{(1-u)^2(1+\chi\int_\Omega\frac{dy}{1-u(y,t)})^2} 
\quad\quad\mbox{in }\Omega\times (0,T)\\
u=&0\qquad\qquad\qquad \qquad\qquad\qquad\quad\quad\quad\quad
\mbox{on }\partial\Omega\times (0,T)\\
u(x,0)=&u_0\qquad\quad\qquad\qquad\qquad\qquad\qquad\quad\quad\,\,
\mbox{ in }\Omega
\end{aligned}\right.
\end{equation*}
in $\Omega\times (0,T)$ if $u\in C^{2,1}(\Omega\times (0,T))
\cap C(\2{\Omega}\times (0,T))$, $0\le u<1$, satisfies 
$$
\frac{\1 u}{\1 t}=\Delta u
+\frac{\lambda f}{(1-u)^2(1+\chi\int_\Omega\frac{dy}{1-u(y,t)})^2} 
\quad\quad\mbox{in }\Omega\times (0,T)
$$
($\le$, $\ge$ respectively) in the classical sense with $u(x,t)=0$ 
($\le$, $\ge$ respectively) on $\1\Omega\times (0,T)$, 
$$
\sup_{\2{\Omega}\times [0,T']}u(x,t)<1\quad\forall 0<T'<T,
$$ 
and 
\begin{equation}
\|u(\cdot,t)-u_0\|_{L^1(\Omega)}\to 0\quad\text{ as }t\to 0.
\end{equation}

Let $\mu_1$ be the first positive eigenvalue and $\phi_1$ be the first 
positive eigenfunction of $-\Delta$ which satisfies $\int_{\Omega}
\phi_1\,dx=1$. For any solution $u$ of ($P_{\lambda}$) we define the 
quenching time $T_{\lambda}>0$ as the time which satisfies 
$$\left\{\aligned
\sup_{\Omega}u(x,t)&<1\quad\forall 0<t<T_{\lambda}\\
\lim_{t\nearrow T_{\lambda}}\sup_{\Omega}u(x,t)&=1.
\endaligned\right.
$$
We say that $u$ has a finite quenching time if $T_{\lambda}<\infty$ and
we say that $u$ quenches at time infinity if $T_{\lambda}=\infty$.

\section{Properties of Steady-states}
\setcounter{equation}{0}
\setcounter{thm}{0}

In this section we will prove the existence of constants $\lambda_N
\ge\lambda_N^{\ast}>0$ such that ($S_{\lambda}^N$) 
has a solution for any $0\le\lambda<\lambda_N^{\ast}$ and ($S_{\lambda}^N$) 
has no solution for any $\lambda>\lambda_N$. For any solution $w$ of 
($S_\lambda$) we let $L_{w,\lambda}=-\Delta-\frac{2\lambda}{(1-w)^3}$ be 
the linearized operator at $w$ and let $\mu_{1,\lambda}(w)$ be 
the first eigenvalue of $L_{w,\lambda}$. We recall a
result of N.~Ghoussoub and Y.~Guo \cite{GG1}.

\begin{thm}(Theorem 1.3 and Theorem 2.1 of \cite{GG1})
There exists a constant $\lambda^{\ast}>0$ such that the following holds.
\begin{enumerate}\rm
\item[(i)] For any $0\le\lambda<\lambda^{\ast}$ there exists a unique minimal 
solution $0\le w_{\lambda}<1$ of ($S_{\lambda}$) such that 
$\mu_{1,\lambda}(w_{\lambda})>0$.
Moreover for each $x\in\Omega$ the function $\lambda\to w_{\lambda}(x)$
is strictly increasing and differentiable on $(0,\lambda^{\ast})$.
\item[(ii)] $\forall\lambda>\lambda^{\ast}$ there is no solution of 
($S_{\lambda}$).
\item[(iii)] Let
\begin{equation}
w_{\ast}=\lim_{\lambda\nearrow\lambda^{\ast}}w_{\lambda}.
\end{equation}
Then $0\le w_{\ast}\le 1$ in $\Omega$. If $1\le n\le 7$, then 
$\sup_{\lambda\in (0,\lambda^{\ast})}
\|w_{\lambda}\|_{\infty}<1$  and $w_{\ast}\in C^{2,\alpha}(\2{\Omega})$ is a 
solution of ($S_{\lambda^{\ast}}$) such that 
$\mu_{1,\lambda^{\ast}}(w_{\ast})=0$.
\end{enumerate} 
\end{thm}

We will now let $\chi>0$, $w_{\lambda}$ be the minimal solution
of ($S_{\lambda}$) given by Theorem 1.1 for any $0<\lambda<\lambda_{\ast}$,
and $w_{\ast}$ be given by (1.1) for the rest of the paper.
Note that by \cite{GPW},
$$ 
\int_{\Omega}\frac{dx}{1-w_{\ast}}<+\infty.
$$

\begin{thm}
Let $\lambda_{\ast}=\lambda^{\ast}(1+\chi\int_{\Omega}
\frac{dx}{1-w_{\ast}})^2$. For any $0\le\lambda<\lambda_{\ast}$ there 
exists a unique constant $\mu_1\in [0,\lambda^{\ast})$ given by 
\begin{equation}
\mu_1\biggl(1+\chi\int_{\Omega}\frac{dy}{1-w_{\mu_1}(y)}\biggr)^2
=\lambda
\end{equation}
such that $w_{\mu_1}$ is a solution of ($S^N_\lambda$). When $1\le n\le 7$, 
the same conclusion holds for any $0\le\lambda\le\lambda_{\ast}$. 
\end{thm}
\begin{proof}
Existence of solution of ($S^N_\lambda$) for $0\le\lambda
<\lambda^{\ast}(1+\chi|\Omega|)^2$ is obtained by 
F.H.~Lin and Y.~Yang in \cite{LY} using a fixed point argument. Here we
will give a simple proof which extends their existence result to the case  
$0\le\lambda<\lambda_{\ast}$. Note 
that $\lambda_{\ast}>\lambda^{\ast}(1+\chi|\Omega|)^2$. When $\lambda=0$, 
the function $v\equiv 0$ is a solution of ($S^N_\lambda$).
Let $0<\lambda<\lambda_{\ast}$ and let
\begin{equation*}
h(\mu)=\mu\left(1+\chi\int_{\Omega}\frac{dx}{1-w_{\mu}}\right)^2
\quad\forall 0\le\mu\le\lambda^{\ast}.
\end{equation*}
Then by Theorem 1.1 $h(\mu)$ is a strictly monotone increasing continuous
function of $\mu\in [0,\lambda^{\ast}]$ and $h(\lambda^{\ast})>\lambda
>h(0)=0$. By the intermediate value theorem there exists a unique
$\mu_1\in (0,\lambda^{\ast})$ satisfying (1.2). Then $w_{\mu_1}$
satisfies ($S^N_\lambda$). When $1\le n\le 7$, by Theorem 1.1 and a 
similar argument the same conclusion holds for any 
$0\le\lambda\le\lambda_{\ast}$.
\end{proof}

Let $\lambda_N^{\ast}=\sup\{\lambda_0>0:(S^N_\lambda)\mbox{ has a solution }
\forall 0\le\lambda\le\lambda_0\}$ and
$$
D_1=\{\lambda_0>0:(S^N_\lambda)\mbox{ has no solution for any }
\lambda\ge\lambda_0\}.
$$
Let $\lambda_N=\inf_{\lambda_0\in D_1}\lambda_0$ if $D_1\ne\phi$ and 
$\lambda_N=\infty$ if $D_1=\phi$. Then by Theorem 1.2, 
$$
\lambda_N\ge\lambda^{\ast}_N
\ge\lambda_{\ast}=\lambda^{\ast}(1+\chi\int_{\Omega}
\frac{dx}{1-w_{\ast}})^2>0.
$$

\begin{prop}
Suppose $\Omega\subset\R^n$ is a strictly convex smooth bounded domain such
that $x\cdot\nu\ge a>0$ for $x\in\1\Omega$ where $\nu$ is the unit outer
normal to $\1\Omega$ at $x$, then $\lambda_N^{\ast}\le\lambda_N<\infty$.
Moreover for any $n\ge 2$,
\begin{equation}
\lambda_N^{\ast}\le\lambda_N\le\frac{(n+2)^2|\1\Omega|}{8an}
\left(\chi(2+\chi|\Omega|)+\frac{1}{|\Omega|}\right).
\end{equation}
\end{prop}
\begin{proof}
Let $\lambda>0$. Suppose $v$ is a solution of ($S^N_{\lambda}$). 
We first claim that there exist constants $C_1>0$ and $\delta>0$ such that
\begin{equation}
\int_{\Omega}\frac{dx}{(1-v)^2}
\le C_1\int_{\Omega\setminus\Omega_{\delta}}\frac{dx}{(1-v)^2}.
\end{equation}
We will use a modification of the proof of Theorem 3.1 of \cite{GW} and 
Theorem 2(a) of \cite{Z} to prove the claim. By the results of 
\cite{GNN}, \cite{GWe}, \cite{T}, there exist constants $\delta>0$ and 
$\alpha>0$
such that $v(x-s\nu_1)$ is an increasing function of $s\in [0,2\delta]$ for 
any $x\in\1\Omega$ and $\nu_1\in\R^n$ satisfying $|\nu_1|=1$ and 
$\nu_1\cdot\nu(x)\ge\alpha$ where $\nu(x)$ is the unit outward normal to
$\1\Omega$ at $x$. Moreover there exists $a_1>0$ such that for any 
$y\in\Omega_{\delta}$ there exists a fixed-sized cone $\Gamma(y)\subset
\Omega_{2\delta}$ with vertex at $y$ such that $|\Gamma(y)\setminus 
\Omega_{\delta}|\ge a_1$ and $v(z)\ge v(y)$ for any $z\in \Gamma(y)$. Then
\begin{align*}
&\frac{1}{(1-v(y))^2}\le\frac{1}{|\Gamma(y)\setminus\Omega_{\delta}|}
\int_{\Gamma(y)\setminus\Omega_{\delta}}\frac{dz}{(1-v(z))^2}
\le\frac{1}{a_1}\int_{\Omega\setminus \Omega_{\delta}}\frac{dz}{(1-v(z))^2}
\quad\forall y\in\Omega_{\delta}\\
\Rightarrow\quad&\int_{\Omega_{\delta}}\frac{dy}{(1-v)^2}
\le\frac{|\Omega|}{a_1}\int_{\Omega\setminus\Omega_{\delta}}
\frac{dx}{(1-v)^2}.
\end{align*}
and (1.4) follows. Multiplying
($S^N_{\lambda}$) by $\phi_1$ and integrating over $\Omega$,
\begin{equation}
\mu_1\ge\mu_1\int_{\Omega}v\phi_1\,dx=-\int_{\Omega}\phi_1\Delta v\,dx
=\lam{\frac{\int_{\Omega}\frac{\phi_1}{(1-v)^2}\,dx}
{(1+\chi\int_{\Omega}\frac{dx}{1-v})^2}}.
\end{equation}
Now by (1.4),
\begin{align}
\left(1+\chi\int_{\Omega}\frac{dx}{1-v}\right)^2
=&1+2\chi\int_{\Omega}\frac{dx}{1-v}+\chi^2\left(
\int_{\Omega}\frac{dx}{1-v}\right)^2\nonumber\\
\le&1+2\chi|\Omega|^{\frac{1}{2}}\left(\int_{\Omega}\frac{dx}{(1-v)^2}
\right)^{\frac{1}{2}}+\chi^2|\Omega|
\int_{\Omega}\frac{dx}{(1-v)^2}\nonumber\\
\le&2\biggl(1+\chi^2|\Omega|
\int_{\Omega}\frac{dx}{(1-v)^2}\biggr)\nonumber\\
\le&2\biggl(1+C_1\chi^2|\Omega|
\int_{\Omega\setminus\Omega_{\delta}}\frac{dx}{(1-v)^2}\biggr)\nonumber\\
\le&2+C_2\chi^2\int_{\Omega}\frac{\phi_1}{(1-v)^2}\,dx
\end{align}
where $C_2=2C_1|\Omega|/\min_{\Omega\setminus\Omega_{\delta}}\phi_1$.
Since $s/(2+C_2\chi^2s)$ is a monotone increasing function of $s\ge 0$, 
by (1.5) and (1.6),
$$
\mu_1\ge\lambda\frac{\int_{\Omega}\frac{\phi_1}{(1-v)^2}\,dx}
{2+C_2\chi^2\int_{\Omega}\frac{\phi_1}{(1-v)^2}\,dx}
\ge\lambda\frac{\int_{\Omega}\phi_1\,dx}
{2+C_2\chi^2\int_{\Omega}\phi_1\,dx}
=\frac{\lambda}{2+C_2\chi^2}.
$$
Hence 
$$
\lambda_N^{\ast}\le\lambda_N\le\mu_1(2+C_2\chi^2)<\infty. 
$$
We will next use a modification of the proof of Proposition 3.3 of \cite{GG1}
to prove (1.3). Suppose $n\ge 2$. By the Pohozaev's identity \cite{N},
\begin{eqnarray}
n\lambda\int_{\Omega}G(v)\,dx-\frac{n-2}{2}\lambda\int_{\Omega}vg(v)\,dx
&=&\frac{1}{2}\int_{\1\Omega}(x\cdot\nu)\left(\frac{\1 v}{\1\nu}\right)^2\,dx
\nonumber\\
&\ge&\frac{a}{2|\1\Omega|}\left(\int_{\1\Omega}\frac{\1 v}{\1\nu}\,dx
\right)^2\nonumber\\
&=&\frac{a}{2|\1\Omega|}\left(-\int_{\1\Omega}\Delta v\,dx\right)^2\nonumber\\
&=&\frac{a\lambda^2}{2b_v^2|\1\Omega|}\left(\int_{\Omega}\frac{dx}{(1-v)^2}
\right)^2
\end{eqnarray}
where $b_v=(1+\chi\int_{\Omega}1/(1-v)\,dx)^2$, $g(v)=1/b_v(1-v)^2$ and 
$$
G(v)=\frac{v}{b_v(1-v)}.
$$
Then by (1.7) and an argument similar to the proof of Proposition 3.3 of 
\cite{GG1},
\begin{eqnarray*}
\frac{a\lambda^2}{|\1\Omega|}\left(\int_{\Omega}\frac{dx}{(1-v)^2}
\right)^2&\le&\lambda b_v\int_{\Omega}\frac{v(n+2-2nv)}{(1-v)^2}
\,dx\\
&\le&\lambda\frac{(n+2)^2}{8n}\left(1+\chi\int_{\Omega}\frac{dx}{1-v}
\right)^2\int_{\Omega}\frac{dx}{(1-v)^2}
\end{eqnarray*}
Hence
\begin{eqnarray*}
\frac{a\lambda}{|\1\Omega|}\int_{\Omega}\frac{dx}{(1-v)^2}
&\le&\frac{(n+2)^2}{8n}\left(1+\chi|\Omega|^{\frac{1}{2}}\left(
\int_{\Omega}\frac{dx}{(1-v)^2}
\right)^{\frac{1}{2}}\right)^2\\
&\le&\frac{(n+2)^2}{8n}\left(1+2\chi|\Omega|^{\frac{1}{2}}
\left(\int_{\Omega}\frac{dx}{(1-v)^2}
\right)^{\frac{1}{2}}+\chi^2|\Omega|\int_{\Omega}\frac{dx}{(1-v)^2}
\right).
\end{eqnarray*}
Thus
\begin{eqnarray*}
\lambda_N^{\ast}
\le\lambda_N&\le&\frac{(n+2)^2|\1\Omega|}{8an}\left(\left(
\int_{\Omega}\frac{dx}{(1-v)^2}\right)^{-1}+2\chi|\Omega|^{\frac{1}{2}}
\left(\int_{\Omega}\frac{dx}{(1-v)^2}\right)^{-\frac{1}{2}}
+\chi^2|\Omega|\right)\\
&\le&\frac{(n+2)^2|\1\Omega|}{8an}\left(\chi(2+\chi|\Omega|)
+\frac{1}{|\Omega|}\right)
\end{eqnarray*}
and (1.3) follows.
\end{proof}

\section{Properties of the nonlocal parabolic MEMS}
\setcounter{equation}{0}

In this section we will prove the local existence and uniqueness of the 
nonlocal parabolic MEMS ($P_{\lambda}$). We also obtain the energy estimates
for the solutions of ($P_{\lambda}$). 

\begin{thm}
Suppose $-b_1\le u_0\in L^1(\Omega)$ satisfies (0.1) for some constants 
$b_1\ge 0$ and $0<a<1$. Let $\lambda>0$, $\chi>0$, and let $u_1$, $u_2$, 
be solutions of ($P_{\lambda}$) in $\Omega\times (0,T)$. Then 
$u_1\equiv u_2$ in $\Omega\times (0,T)$.
\end{thm}
\begin{proof}
We will use a modification of the technique of Dahlberg and Kenig \cite{DK}
and K.M.~Hui \cite{H} to prove the lemma. By reducing $T$ slightly we may 
assume without loss of generality that
\begin{equation}
\sup_{\Omega\times (0,T)}u_i\le b_2<1\quad\forall i=1,2,
\end{equation}
for some constant $0<b_2<1$. Note that since $u_i$ is a 
supersolution of the heat equation in $\Omega\times (0,T)$ for $i=1,2$,
by the maximum principle,
\begin{equation}
u_i\ge -b_1\quad\mbox{ in }\Omega\times (0,T)\quad\forall i=1,2.
\end{equation}
Let $h\in C_0^{\infty}(\Omega)$ be such that $0\le h\le 1$. For any 
$t_1\in (0,T)$, let $\eta$ be the solution of 
\begin{equation}
\left\{\begin{aligned}
\eta_t+\Delta\eta+H\eta&=0\qquad\text{ in }\Omega\times (0,t_1)\\
\eta&=0\qquad\text{ on }\1\Omega\times (0,t_1)\\
\eta (x,t_1)&=h(x)\quad\text{ in }\Omega
\end{aligned}\right.
\end{equation}  
where
\begin{equation}
H(x,t)
=\left\{\begin{aligned}
&\lambda\left(\frac{(1-u_1)^{-2}-(1-u_2)^{-2}}{u_1-u_2}\right)
\frac{1}{(1+\chi\int_\Omega\frac{dy}{1-u_1(y,t)})^2}
\quad\text{ if }u_1(x,t)\ne u_2(x,t)\\
&\frac{2\lambda}{(1-u_1)^3}
\frac{1}{(1+\chi\int_\Omega\frac{dy}{1-u_1(y,t)})^2}\qquad
\qquad\qquad\qquad\qquad\text{if }u_1(x,t)=u_2(x,t).\end{aligned}\right.
\end{equation} 
By the maximum principle $\eta\ge 0$. Then
\begin{align}
&\int_{\Omega}(u_1-u_2)(x,t_1)h(x)\,dx\nonumber\\
=&\int_0^{t_1}\int_{\Omega}\frac{\1}{\1 t}[(u_1-u_2)\eta]\,dx\,dt\nonumber\\
=&\int_0^{t_1}\int_{\Omega}[(u_1-u_2)_t\eta+(u_1-u_2)\eta_t]\,dx\,dt
\nonumber\\
=&\int_0^{t_1}\int_{\Omega}\biggl[\eta\Delta (u_1-u_2)+(u_1-u_2)\eta_t
\nonumber\\
&\qquad+\lambda\eta\biggl(\frac{1}{(1-u_1)^2(1+\chi\int_\Omega
\frac{dy}{1-u_1(y,t)})^2}-\frac{1}{(1-u_2)^2(1+\chi\int_\Omega
\frac{dy}{1-u_2(y,t)})^2}\biggr)\biggr]\,dx\,dt\nonumber\\
=&\int_0^{t_1}\int_{\Omega}(u_1-u_2)[\eta_t+\Delta\eta+H\eta]\,dx\,dt
\nonumber\\
&\qquad+\int_0^{t_1}\int_{\Omega}\frac{\lambda\eta}{(1-u_2)^2}
\biggl(\frac{1}{(1+\chi\int_\Omega\frac{dy}{1-u_1(y,t)})^2}
-\frac{1}{(1+\chi\int_\Omega\frac{dy}{1-u_2(y,t)})^2}\biggr)\,dx\,dt
\end{align}
Hence by (2.1), (2.2), (2.3) and (2.5),
\begin{align}
&\int_{\Omega}(u_1-u_2)(x,t_1)h(x)\,dx\nonumber\\
=&\int_0^{t_1}\int_{\Omega}\frac{\lambda\chi\eta}{(1-u_2)^2}\cdot
\frac{\int_{\Omega}\frac{u_2-u_1}{(1-u_2)(1-u_1)}\,dy\left(2+
\chi\int_\Omega\frac{dy}{1-u_1(y,t)}
+\chi\int_\Omega\frac{dy}{1-u_2(y,t)}\right)}
{(1+\chi\int_\Omega\frac{dy}{1-u_1(y,t)})^2
(1+\chi\int_\Omega\frac{dy}{1-u_2(y,t)})^2}\,dxdt\nonumber\\
\le&C_1\|\eta\|_{\infty}\int_0^{t_1}\int_\Omega(u_2-u_1)_+(y,t)\,
dydt
\end{align}
for some constant $C_1>0$ depending on $b_1$, $b_2$, $\lambda$ and $\chi$. 
By (2.1), (2.2) and (2.4) there exists a constant $C_2>0$ such that
$\sup_{\Omega\times (0,T)}|H(x,t)|\le C_2$. Then
\begin{align}
&\eta_t+\Delta\eta+C_2\eta\ge 0\qquad\quad\text{in }\Omega\times (0,t_1)\\
&(e^{C_2t}\eta)_t+\Delta (e^{C_2t}\eta)\ge 0\quad\text{ in }\Omega\times 
(0,t_1).
\end{align}
Hence by the maximum principle,
\begin{equation}
\eta (x,t)\le e^{C_2t}\eta (x,t)\le\max_{\2{\Omega}} (e^{C_2t_1}
\eta (x,t_1))=e^{C_2t_1}\|h\|_{L^{\infty}}\le e^{C_2T}\quad\forall
x\in\Omega,0<t<t_1.
\end{equation}
By (2.6) and (2.9),
\begin{equation}
\int_{\Omega}(u_1-u_2)(x,t_1)h(x)\,dx
\le C_1e^{C_2T}\int_0^{t_1}\int_\Omega(u_2-u_1)_+(y,t)dy\,dt
\end{equation}
We now choose a sequence of smooth function $\{h_i\}$, $0\le h_i\le 1$,
such that $h_i(x)$ converges a.e. to the characteristic function of the set
$\{x:u_1(x,t_1)>u_2(x,t_1)\}$ as $i\to\infty$. Letting $i\to\infty$ in 
(2.10) we get
\begin{equation}
\int_{\Omega}(u_1-u_2)_+(x,t_1)\,dx
\le C_1e^{C_2T}\int_0^{t_1}\int_\Omega(u_2-u_1)_+(y,t)\,dy\,dt
\quad\forall 0<t_1<T.
\end{equation}
Interchanging the role of $u_1$ and $u_2$,
\begin{equation}
\int_{\Omega}(u_2-u_1)_+(x,t_1)\,dx
\le C_1e^{C_2T}\int_0^{t_1}\int_\Omega(u_1-u_2)_+(y,t)\,dy\,dt
\quad\forall 0<t_1<T.
\end{equation}
By (2.11) and (2.12),
\begin{equation}
\int_{\Omega}|u_1-u_2|(x,t_1)\,dx
\le C_1e^{C_2T}\int_0^{t_1}\int_\Omega|u_1-u_2|(y,t)\,dy\,dt
\quad\forall 0<t_1<T.
\end{equation}
Let
$$
y(t)=\int_0^t\int_\Omega|u_1-u_2|(z,s)\,dz\,ds.
$$
Then by (2.13),
\begin{align*}
&y'(t)\le C_1e^{C_2T}y(t)\quad\forall 0<t<T\\
\Rightarrow\quad&y(t)\le e^{C_1e^{C_2T}T}y(0)=0\quad\forall 0<t<T\\
\Rightarrow\quad&y(t)\equiv 0\quad\forall 0<t<T
\end{align*}
and the theorem follows.
\end{proof}

\begin{thm}
Let $-b\le u_0\in L^1(\Omega)$ satisfies (0.1) for some constants 
$b\ge 0$ and $0<a<1$. Then for any $\lambda>0$ and $\chi>0$ there 
exists $T>0$ such that ($P_{\lambda}$) has a solution $u\ge -b$
in $\Omega\times (0,T)$ which satisfies
\begin{equation}
u(x,t)=\int_{\Omega}G(x,y,t)u_0(y)\,dy+\lambda\int_0^t\int_{\Omega}
\frac{G(x,y,t-s)}{(1-u(y,s))^2(1+\chi\int_{\Omega}\frac{dz}{1-u(z,s)})^2}
\,dy\,ds
\end{equation}
for all $(x,t)\in\Omega\times (0,T)$ where $G(x,y,t)$ is the Dirichlet
Green function for the heat equation in $\Omega\times (0,T)$.
\end{thm}
\begin{proof}
We will use a modification of the proof of Theorem 2.5 of \cite{H} to prove 
the theorem. 

\noindent $\underline{\text{\bf Case 1}}$: $-b\le u_0\in 
C_0^{\infty}(\Omega)$ satisfies (0.1).

Let 
\begin{equation}
T=\frac{(1-a)^3}{16\lambda},
\end{equation}
\begin{equation}
u_1(x,t)=\int_{\Omega}G(x,y,t)u_0(y)\,dy+\lambda\int_0^t\int_{\Omega}
\frac{G(x,y,t-s)}{(1-u_0(y))^2(1+\chi\int_{\Omega}\frac{dz}{1-u_0(z)})^2}
\,dy\,ds
\end{equation}
and
\begin{equation}
u_{k+1}(x,t)=\int_{\Omega}G(x,y,t)u_0(y)\,dy+\lambda\int_0^t\int_{\Omega}
\frac{G(x,y,t-s)}{(1-u_k(y,s))^2(1+\chi\int_{\Omega}\frac{dz}{1-u_k(z,s)})^2}
\,dy\,ds
\end{equation}
for all $x\in\2{\Omega}, 0<t<T,k\ge 1$.
Let
\begin{equation}
\4{u}_1(x,t)=\int_{\Omega}G(x,y,t)u_0(y)\,dy+\lambda\int_0^t\int_{\Omega}
\frac{G(x,y,t-s)}{(1-u_0(y))^2}\,dy\,ds
\end{equation}
and
\begin{equation}
\4{u}_{k+1}(x,t)=\int_{\Omega}G(x,y,t)u_0(y)\,dy+\lambda\int_0^t\int_{\Omega}
\frac{G(x,y,t-s)}{(1-\4{u}_k(y,s))^2}\,dy\,ds
\end{equation}
for all $x\in\2{\Omega}, 0<t<T,k\ge 1$. Then 
$$
u_1\le\4{u}_1\quad\forall x\in\2{\Omega}, 0\le t\le T.
$$
Suppose
\begin{equation}
u_k\le\4{u}_k\quad\forall x\in\2{\Omega}, 0\le t\le T,
\end{equation}
holds for some $k\ge 1$. Then by (2.17), (2.19) and (2.20),
\begin{equation}
u_{k+1}\le\4{u}_{k+1}\quad\forall x\in\2{\Omega}, 0\le t\le T.
\end{equation}
Hence by induction (2.20) holds for all $k\ge 1$. Since by the
proof of Theorem 2.5 of \cite{H},
$$
\4{u}_k\le\frac{1+a}{2}\quad\forall x\in\2{\Omega}, 0\le t\le T,k\ge 1,
$$
by (2.20),
\begin{equation}
u_k\le\frac{1+a}{2}\quad\forall x\in\2{\Omega}, 0\le t\le T,k\ge 1.
\end{equation}
Let
\begin{equation}
q(x,t)=\int_{\Omega}G(x,y,t)u_0(y)\,dy
\end{equation}
Then $q$ is the solution of the problem
\begin{equation}\left\{\aligned 
\1_tq&=\Delta q\qquad\text{ in }\Omega\times (0,\infty)\\
q&=0\qquad\quad\,\mbox{on }\1\Omega\times (0,\infty)\\
q(x,0)&=u_0(x)\quad\text{ in }\Omega.\endaligned\right.
\end{equation}
By (2.16), (2.17) and (2.23),
\begin{equation}
u_k(x,t)\ge q(x,t)\ge -b\quad\forall x\in\2{\Omega}, 0\le t\le T, k\in\Z^+.
\end{equation}
Since 
\begin{equation}
q(\cdot,t)\to u_0\quad\text{ in }L^1(\Omega)\quad\text{ as }t\to 0,
\end{equation}
by (2.16), (2.17), (2.22) and (2.25),
\begin{equation}
u_k(\cdot,t)\to u_0\quad\text{ in }L^1(\Omega)\quad\text{ as }t\to 0.
\end{equation}
By (2.16), $u_1$ is continuously differentiable in $x$ and $t$. Then by
(2.16), (2.17), (2.22), (2.25) and standard parabolic theory \cite{LSU}, 
\cite{F}, $u_k\in C^{2,1}
(\2{\Omega}\times (0,T])$ for all $k\ge 2$. Then by (2.16), (2.17) and 
(2.22), $\forall k\ge 2$, $u_k$ satisfies
\begin{equation}\left\{\aligned
\frac{\1 u_k}{\1 t}=&\Delta u_k
+\frac{\lambda}{(1-u_{k-1})^2(1+\chi\int_{\Omega}\frac{dy}{1-u_{k-1}(y,t)})^2}
\quad\text{ in }\Omega\times (0,T)\\
u_k(x,t)&=0\qquad\qquad\qquad\qquad\qquad\qquad\qquad
\,\,\text{ on }\1\Omega\times (0,T)\\
u_k(x,0)&=u_0(x)\quad\qquad\qquad\qquad\qquad\qquad\qquad\text{in }\Omega.
\endaligned\right.
\end{equation}
By (2.22), (2.25), (2.28) and the parabolic Schauder estimates \cite{LSU}, 
the sequence $\{u_k\}_{k=2}^{\infty}$ are uniformly Holder continuous 
on $\2{\Omega}\times [0,T]$. Then by (2.22), (2.25), (2.28) and
the Schauder estimates 
(\cite{LSU},\cite{F}) $\{u_k\}_{k=2}^{\infty}$ are uniformly bounded in
$C^{2,1}(K)$ for any compact subset $K\subset\2{\Omega}\times (0,T]$. 
By the Ascoli-Arzel\'a theorem and a 
diagonalization argument $\{u_k\}_{k=2}^{\infty}$ has a subsequence
which we may assume without loss of generality to be the sequence itself
which converges uniformly in $C^{2,1}(K)$ to some 
function $u$ for any
compact subset $K\subset\2{\Omega}\times (0,T]$ as $k\to\infty$. Then by 
(2.16), (2.17), (2.22) and (2.25) $u$ satisfies (2.14),
\begin{equation}
-b\le q(x,t)\le u(x,t)\le\frac{1+a}{2}\quad\forall x\in\2{\Omega}, 0<t\le T,
\end{equation}
and
$$\left\{\aligned
\frac{\1 u}{\1 t}=&\Delta u
+\frac{\lambda}{(1-u)^2(1+\chi\int_{\Omega}\frac{dy}{1-u(y,t)})^2}
\quad\text{ in }\Omega\times (0,T)\\
u(x,t)&=0\qquad\qquad\qquad\qquad\qquad\qquad
\quad\text{on }\1\Omega\times (0,T)\\
\endaligned\right.
$$  
By (2.14), (2.26) and (2.29), $u$ satisfies (0.2).
Hence $u$ is a solution of ($P_{\lambda}$) in $\Omega\times (0,T)$.
 
\noindent $\underline{\text{\bf Case 2}}$: $-b\le u_0\in L^1(\Omega)$
satisfies (0.1). 

We choose a sequence of function $\{u_{0,k}\}_{k=1}^{\infty}\in 
C_0^{\infty}(\Omega)$ such that $-b\le u_{0,k}\le a$ in $\Omega$ for all 
$k\ge 1$ and $u_{0,k}$ converges to $u_0$ in $L^1(\Omega)$ and a.e. as 
$k\to\infty$. For any $k\in\mathbb{Z}^+$, by case 1 there exists a 
solution $u_k$ of ($P_{\lambda}$) in $\Omega\times (0,T)$ with initial 
value $u_{0,k}$ which satisfies (2.22), (2.27) with $u_0$ there being 
replaced by $u_{0,k}$,  
$$
u_k(x,t)\ge\int_{\Omega}G(x,y,t)u_{0,k}(y)\,dy\ge -b\quad\forall
x\in\2{\Omega},0\le t\le T,k\in\Z^+, 
$$
and
\begin{equation}
u_k(x,t)=\int_{\Omega}G(x,y,t)u_{0,k}(y)\,dy
+\lambda\int_0^t\int_{\Omega}\frac{G(x,y,t-s)}{(1-u_k(y,s))^2
(1+\chi\int_{\Omega}\frac{dz}{1-u_k(z,s)})^2}\,dy\,ds
\end{equation}
for any $x\in\2{\Omega},0<t<T$ where $T$ is given by (2.15). Then by an 
argument similar to \cite{H},
the sequence $\{u_k\}_{k=1}^{\infty}$ are uniformly bounded in 
$C^{2,1}(K)$ for any compact subset $K\subset\2{\Omega}\times (0,T]$. 
Moreover $\{u_k\}_{k=1}^{\infty}$ has a subsequence
which we may assume without loss of generality to be the sequence itself
which converges uniformly in $C^{2,1}(K)$ to a solution 
$u$ of ($P_{\lambda}$) which satisfies (2.14) and (2.29) with $q(x,t)$ 
being given by (2.23) and the theorem follows.
\end{proof}

\begin{thm}
Let $\lambda>0$, $\chi>0$, and let $0\le u_0\in L^1(\Omega)$ satisfy 
(0.1) for some constant $0<a<1$. Let $u$ be a global solution 
of ($P_{\lambda}$). Then $u$ satisfies
\begin{equation}
\int_{t_0}^T\int_{\Omega}u_t^2\,dxdt
+\frac{1}{2}\int_{\Omega}|\nabla u(x,T)|^2
\,dx\le\frac{1}{2}\int_{\Omega}|\nabla u(x,t_0)|^2\,dx
+\frac{\lambda}{\chi(1+\chi|\Omega|)}
\end{equation}
for any $T>t_0>0$. If $u_0\equiv 0$ on $\Omega$, then 
\begin{equation}
\int_0^T\int_{\Omega}u_t^2\,dxdt+\frac{1}{2}\int_{\Omega}
|\nabla u(x,T)|^2\,dx\le\frac{\lambda}{\chi(1+\chi|\Omega|)}
\quad\forall T>0.
\end{equation}
\end{thm}
\begin{proof} 
By Theorem 2.1 and Theorem 2.2 $u$ is uniquely given by (2.14) in 
$\Omega\times (0,\infty)$. Hence $u\ge 0$ in $\Omega\times (0,\infty)$.
Then
\begin{align*}
\int_{t_0}^t\int_{\Omega}u_t^2\,dx\,dt=&\int_{t_0}^t
\int_{\Omega}u_t\Delta u\,dx\,dt
+\lambda\int_{t_0}^t\int_{\Omega}\frac{u_t}{(1-u)^2(1+\chi\int_{\Omega}
\frac{dy}{1-u})^2}\,dx\,dt\\
=&-\frac{1}{2}\int_{t_0}^t\frac{\1}{\1 t}\left(\int_{\Omega}|\nabla u|^2\,dx
\right)dt+\lambda\int_{t_0}^t\frac{\frac{\1}{\1 t}\left(\int_{\Omega}
\frac{dx}{1-u}\right)}{(1+\chi\int_{\Omega}
\frac{dy}{1-u})^2}\,dt\\
=&\frac{1}{2}\int_{\Omega}|\nabla u(x,t_0)|^2\,dx-\frac{1}{2}\int_{\Omega}
|\nabla u(x,t)|^2\,dx+\frac{\lambda}{\chi(1+\chi\int_{\Omega}
\frac{dy}{1-u(y,t_0)})}\\
&\qquad-\frac{\lambda}{\chi(1+\chi\int_{\Omega}
\frac{dy}{1-u(y,t)})}\qquad\qquad\forall t\ge t_0>0.
\end{align*}
Hence 
\begin{align}
&\int_{t_0}^t\int_{\Omega}u_t^2\,dx\,dt
+\frac{1}{2}\int_{\Omega}|\nabla u(x,t)|^2
\,dx+\frac{\lambda}{\chi(1+\chi\int_{\Omega}
\frac{dy}{1-u(y,t)})}\nonumber\\
=&\frac{1}{2}\int_{\Omega}|\nabla u(x,t_0)|^2\,dx
+\frac{\lambda}{\chi(1+\chi\int_{\Omega}
\frac{dy}{1-u(y,t_0)})}\nonumber\\
\le&\frac{1}{2}\int_{\Omega}|\nabla u(x,t_0)|^2\,dx
+\frac{\lambda}{\chi(1+\chi|\Omega|)}\quad\forall t\ge t_0>0.\nonumber
\end{align}
and (2.31) follows. If $u_0\equiv 0$ on 
$\Omega$, then by (2.14) and \cite{LSU},
\begin{align}
u(x,t)=&\lambda\int_0^t\int_{\Omega}
\frac{G(x,y,t-s)}{(1-u(y,s))^2(1+\chi\int_{\Omega}\frac{dz}{1-u(z,s)})^2}
\,dy\,ds\nonumber\\
\Rightarrow\quad|\nabla u(x,t)|\le&\lambda\int_0^t\int_{\Omega}
\frac{|\nabla_xG(x,y,t-s)|}{(1-u(y,s))^2
(1+\chi\int_{\Omega}\frac{dz}{1-u(z,s)})^2}
\,dy\,ds\nonumber\\
\le&C\int_0^t\int_{\Omega}|\nabla G(x,y,t-s)|\,dy\,ds\nonumber\\
\le&C\int_0^t\int_{\R^n}\frac{|x-y|}{(t-s)^{\frac{n+1}{2}}}
e^{-c\frac{|x-y|^2}{t-s}}\,dy\,ds\nonumber\\
\le&C t\quad\forall x\in\Omega, 0<t\le 1.
\end{align}
where $C>0$ is a generic constant that changes from line to line.
Letting $t_0\to 0$ in (2.31) by (2.33) we get (2.32) and the theorem follows. 
\end{proof}

\begin{cor}
Let $\lambda>0$ and $\chi>0$. Let $0\le u_0\in L^1(\Omega)$ satisfy (0.1)
for some constant $0<a<1$. Suppose $u$ is a global solution of 
($P_{\lambda}$) which satisfies $0\le u\le b$ in $\Omega\times (0,\infty)$ 
for some constant $0<b<1$. Let $\{t_i\}_{i=1}^{\infty}$ be a sequence 
such that $t_i>t_0$ for all $i\ge 1$ and $t_i\to\infty$ as $i\to\infty$ 
for some constant $t_0>0$. Then the sequence $\{t_i\}_{i=1}^{\infty}$ 
has a subsequence $\{t_i'\}_{i=1}^{\infty}$ such that $u(x,t_i')$ 
converges uniformly in $C^2(\2{\Omega})$ to a solution $v$ of 
($S^N_{\lambda}$) as $i\to\infty$,
\begin{equation}
\int_{\Omega}u_t^2(x,t_i')\,dx\to 0\quad\mbox{ as }i\to\infty  
\end{equation} 
and $v$ satisfies
\begin{equation}
\int_{\Omega}|\nabla v|^2\,dx\le\int_{\Omega}|\nabla u(x,t_0)|^2\,dx
+\frac{2\lambda}{\chi(1+\chi|\Omega|)}.
\end{equation}
If $u_0\equiv 0$ on $\Omega$, then the same conclusion holds with
$t_0=0$.
\end{cor}
\begin{proof}
Since $0\le u\le b$ in $\Omega\times (0,\infty)$, by the parabolic 
Schauder estimates \cite{LSU} $u\in C^{2+\beta,1+(\beta/2)}(\2{\Omega}
\times [1,\infty))$.  Then by the Ascoli-Arzel\'a theorem there exists a subsequence 
$\{t_i'\}_{i=1}^{\infty}$ of $\{t_i\}_{i=1}^{\infty}$ such that 
$u(x,t_i'+t)$ converges uniformly in $C^2(\2{\Omega})$ to some function 
$v_1$ in $C^{2+\beta,1+(\beta/2)}(\2{\Omega}\times [0,1])$ as $i\to\infty$. 
Let $v(x)=v_1(x,0)$. By Theorem 2.3 (2.31) holds. Suppose there exists a 
constant $\3>0$ and a subsequence $\{t_i''\}_{i=1}^{\infty}$ of 
$\{t_i'\}_{i=1}^{\infty}$ such that
\begin{equation}
\int_{\Omega}u_t^2(x,t_i'')\,dx\ge\3\quad\forall i\in\Z^+.
\end{equation} 
Without loss of generality we may assume that $t_i''\ge t_{i-1}''+2$ for
all $ i\in\Z^+$ and $t_0\ge 2$. Since $u\in C^{2+\beta,1+(\beta/2)}(\2{\Omega}
\times [1,\infty))$, there exists a constant $0<\delta<1$ such that
\begin{equation}
\biggl|\int_{\Omega}u_t^2(x,t)\,dx-\int_{\Omega}u_t^2(x,t')\,dx\biggr|
\le\frac{\3}{2}\quad\forall t,t'\ge t_0,|t-t'|\le\delta.
\end{equation} 
By (2.36) and (2.37),
\begin{align}
&\int_{\Omega}u_t^2(x,t)\,dx\ge\frac{\3}{2}\quad\forall 
|t-t_i''|\le\delta,i\in\Z^+\nonumber\\
\Rightarrow\quad&\int_{t_0}^{\infty}\int_{\Omega}u_t^2(x,t)\,dx\,dt
\ge\sum_{i=1}^{\infty}\int_{t_i''-\delta}^{t_i''+\delta}
\int_{\Omega}u_t^2(x,t)\,dx\,dt=\sum_{i=1}^{\infty}\delta\3
=\infty.\nonumber
\end{align} 
This contradicts (2.31). Hence (2.34) holds.
By (2.31) and an 
argument similar to the proof of Theorem 3.1 of \cite{H} $v_1(x,t)=v(x)$
for all $0\le t\le 1$. Integrating ($P_{\lambda}$) over $(t_i',t_i'+1)$,
\begin{align*}
&u(x,t_i'+1)-u(x,t_i')=\int_{t_i'}^{t_i'+1}\Delta u
+\int_{t_i'}^{t_i'+1}
\frac{\lambda}{(1-u)^2(1+\chi\int_\Omega\frac{dy}{1-u(y,t)})^2}\,dt\\
\Rightarrow\quad&0=\Delta v
+\frac{\lambda}{(1-v)^2(1+\chi\int_\Omega\frac{dy}{1-v(y)})^2}
\quad\mbox{ as }i\to\infty.
\end{align*}
Putting $t=t_i'$ in (2.31) and letting $i\to\infty$ (2.35) follows. If 
$u_0\equiv 0$ on $\Omega$, then (2.32) holds. Putting $t=t_i'$ in (2.32) 
and letting $i\to\infty$ (2.35) holds with $t_0=0$ and the corollary follows.
\end{proof}

\begin{cor}
Let $\lambda>\lambda_N$ and $\chi>0$. Let $u$ be a global solution of 
($P_{\lambda}$). Then either $T_{\lambda}<\infty$ or $u$ quenches at 
time infinity.
\end{cor}

\section{Global existence and asymptotic behaviour of solutions of 
($P_{\lambda}$)}
\setcounter{equation}{0}

In this section we will prove the global existence and asymptotic behaviour 
of solutions of ($P_{\lambda}$) under various boundedness conditions on 
$\lambda$.

\begin{thm}
Let $n=1$, $b>0$, $\Omega=(-b,b)$, $\chi>0$ and $0<\lambda
<\frac{\chi (1+\chi |\Omega|)}{2|\Omega|}$. Then there exists a unique
global solution $u\ge 0$ for ($P_{\lambda})$ in $\Omega\times (0,\infty)$ 
with $u_0=0$ and there exists a solution $v$ of ($S^N_{\lambda}$) which 
satisfies
\begin{equation}
\int_{-b}^bv_x^2\,dx\le\frac{2\lambda}{\chi(1+\chi|\Omega|)}.
\end{equation}
Hence $\lambda_N^{\ast}\ge\chi (1+\chi |\Omega|)/2|\Omega|$.
\end{thm}
\begin{proof}
By Theorem 2.1 we only need to prove existence of global solution of 
($P_{\lambda})$.
By Theorem 2.2 there exists $T'>0$ such that ($P_{\lambda})$ has a solution
$u\ge 0$ in $\Omega\times (0,T')$ with $u_0=0$. Let $T>0$ be the maximal 
time of existence of a solution $u\ge 0$ of ($P_{\lambda})$ in 
$\Omega\times (0,T)$ with $u_0=0$. Suppose $T<\infty$. By Theorem 2.3,
\begin{equation*}
\int_{-b}^bu_x^2(x,t)\,dx\le\frac{2\lambda}{\chi(1+\chi|\Omega|)}
\quad\forall 0\le t<T.
\end{equation*}
Hence $\forall |x|\le b, 0\le t<T$,
$$
u(x,t)=\int_{-b}^xu_x(y,t)\,dy
\le|\Omega|^{\frac{1}{2}}\left(\int_{-b}^bu_x^2(y,t)\,dy\right)^{\frac{1}{2}}
\le\sqrt{\frac{2\lambda|\Omega|}{\chi(1+\chi|\Omega|)}}<1
$$
Let $a_1=\sqrt{\frac{2\lambda|\Omega|}{\chi(1+\chi|\Omega|)}}$, 
$T_1=(1-a_1)^3/16\lambda$, and $T_2=T-\min (T_1/2,T/2)$. By the 
proof of Theorem 2.2 there exists a solution $u_1$ of ($P_{\lambda})$
in $\Omega\times (0,T_1)$ with $u_0=u(x,T_2)$. We then extend $u$ to a
function on $\2{\Omega}\times (0,T_2+T_1)$ be setting $u(x,t)=u_1(x,t-T_2)$
for all $x\in\2\Omega$ and $T_2\le t\le T_2+T_1$. Then $u$ is a solution
of ($P_{\lambda})$ in $\Omega\times (0,T_2+T_1)$ with $u_0=0$. Since 
$T_1+T_2>T$, this contradicts the maximality of $T$. Hence $T=\infty$ and
$u$ is a global solution of ($P_{\lambda})$ with $u_0=0$. By Corollary 2.4
there exists a sequence $t_i\to\infty$ as $i\to\infty$ such that $u(x,t_i)$
converges uniformly on $\2{\Omega}$ to a solution $v$ of ($S^N_{\lambda}$)
which satisfies (2.35) with $t_0=0$ as $i\to\infty$
and the theorem follows.
\end{proof}

\begin{cor}
Let $n=1$, $b>0$ and $\Omega=(-b,b)$. For any  $\chi>0$ and $0<\lambda
<\frac{\chi (1+\chi |\Omega|)}{2|\Omega|}$, let $v_{\lambda,\chi}$ be the
solution of ($S^N_{\lambda}$) constructed in Theorem 3.1. Then
$v_{\lambda,\chi}$ converges uniformly to $0$ on $\2{\Omega}$ as 
$\chi\to\infty$ or $\lambda\to 0$.
\end{cor}
\begin{proof}
By Theorem 3.1 $v_{\lambda,\chi}$ satisfies (3.1). Hence
$$
|v_{\lambda,\chi}(x)|\le\left|\int_{-b}^x(v_{\lambda,\chi})_x\,dy\right|
\le\sqrt{2b}\left(\int_{-b}^b(v_{\lambda,\chi})_x^2\,dy\right)^{\frac{1}{2}}
\le 2\sqrt{\frac{b\lambda}{\chi (1+\chi|\Omega|)}}.
$$
Since the right hand side tends to $0$ uniformly on $\2{\Omega}$ as 
$\chi\to\infty$ or $\lambda\to 0$, the corollary follows.
\end{proof}

We next recall a result of \cite{H}.

\begin{lem}(cf. Theorem 2.1 of \cite{H})
Let $u_{0,1}, u_{0,2}\in L^1(\Omega)$ be such that $0\le u_{0,1}
\le u_{0,2}\le a$ in $\Omega$ for some constant $0<a<1$. Let 
$0\le f\in C(\2{\Omega}\times (0,T))\cap L^{\infty}(\2{\Omega}\times 
(0,T))$. Suppose $u_1$, $u_2$, are nonnegative subsolution and 
supersolution of 
\begin{equation*}\left\{\begin{aligned}
\frac{\1 u}{\1 t}=&\Delta u
+\frac{f(x,t)}{(1-u)^2} 
\quad\quad\mbox{in }\Omega\times (0,T)\\
u=&0\qquad\qquad\qquad\quad\,\,\mbox{ on }\partial\Omega\times (0,T)\\
u(x,0)=&u_0\qquad\qquad\qquad\quad\mbox{ in }\Omega
\end{aligned}\right.
\end{equation*}
in $\Omega\times (0,T)$ with 
initial value $u_0=u_{0,1}, u_{0,2}$, 
respectively which satisfy (2.1) for some constant $0<b_2<1$. Then 
$u_1\le u_2$ in $\2{\Omega}\times (0,T)$.
\end{lem}

\begin{thm}
Let $0<\lambda<\lambda^{\ast}(1+\chi|\Omega|)^2$, $\chi>0$ and 
$\mu=\lambda/(1+\chi|\Omega|)^2$. Let $\mu_0\in [\mu,\lambda^{\ast})$ and
let $u_0\in L^1(\Omega)$ be such that $0\le u_0\le w_{\mu_0}$ in $\Omega$, 
there exists a unique global solution $u$ of ($P_{\lambda}$) such that 
$0\le u\le w_{\mu_0}$ in $\Omega\times (0,\infty)$. If $1\le n\le 7$, 
the same result remains valid for $0<\lambda
\le\lambda^{\ast}(1+\chi|\Omega|)^2$ and $\mu_0=\lambda^{\ast}$.
\end{thm}
\begin{proof}
As before by Theorem 2.1 we only need to prove existence of global solution 
of ($P_{\lambda}$). Let $0<\lambda<\lambda^{\ast}(1+\chi|\Omega|)^2$.
Note that since the inequality
$$
u(x,t)\le |\Omega|^{\frac{1}{2}}\|u_x\|_{L^2(\Omega)}
$$
for function $u$ vanishing at $\1\Omega\times (0,T)$
is only valid when $\Omega$ is a bounded interval in $\R$ and $n=1$,
the argument of the proof of Theorem 3.1 cannot be used here. We will
use another method to prove this theorem.

By Theorem 2.2 there exists $T'>0$ such that ($P_{\lambda}$) has a 
non-negative solution $u$ in $\Omega\times (0,T')$. Let $T>0$ be the 
maximal time of existence of the solution $u$. 
Then $u$ satisfies
$$
\frac{\1 u}{\1 t}\le\Delta u+\frac{\mu_0}{(1-u)^2} 
\quad\quad\mbox{in }\Omega\times (0,T).
$$
Hence by Lemma 3.3 $u\le w_{\mu_0}\le\|w_{\mu_0}\|_{\infty}<1$ on 
$\2{\Omega}\times (0,T)$. Suppose $T<\infty$. Then by the parabolic 
Schauder estimates \cite{LSU} $u\in C^{2,1}(\2{\Omega} 
\times (T/2,T))$. Hence $u$ can be
extended to a function on $\2{\Omega} \times (T/2,T]$ and $u\in 
C^{2,1}(\2{\Omega} \times (T/2,T])$.
By Theorem 2.2 there exists $\delta>0$ such that there exists a solution
$\2{u}$ of ($P_{\lambda}$) in $\Omega\times (0,\delta)$ with initial value
$\2{u}(x,0)=u(x,T)$. Let $u(x,t)=\2{u}(x,t-T)$ for any $x\in\2{\Omega}$
and $T\le t\le T+\delta$. Then $u$ is a solution of ($P_{\lambda}$) in 
$\Omega\times (0,T+\delta)$. This contradicts the maximality of $T$. Hence 
$T=\infty$ and $u$ is a global solution of ($P_{\lambda}$). 

If $1\le n\le 7$, by Theorem 1.1 $w_{\ast}<1$ on $\2{\Omega}$. Then
the same argument as before also works for the case
$0<\lambda\le\lambda^{\ast}(1+\chi|\Omega|)^2$ and $\mu_0=\lambda^{\ast}$
and the theorem follows.
\end{proof}

\begin{thm}
Let $\Omega\subset\R^n$ be a smooth bounded domain. Let $0<\lambda
<\lambda^{\ast}(1+\chi|\Omega|)^2$, $\chi>0$ and $\mu
=\lambda/(1+\chi|\Omega|)^2$. Let $\mu_0\in [\mu,\lambda^{\ast})$ and
let $u_0\in L^1(\Omega)$ be such that $0\le u_0\le w_{\mu_0}$ in $\Omega$. 
Suppose $u$ is the unique global 
solution of ($P_{\lambda}$) given by Theorem 3.4 which satisfies 
$0\le u\le w_{\mu_0}$ in $\Omega\times (0,\infty)$. Then there exists a 
constant $0<\mu_1\le\mu_0$ given uniquely by (1.2) such that 
$u(\cdot,t)$ converges uniformly in $C^2(\2{\Omega})$ to $w_{\mu_1}$ 
as $t\to\infty$. If $1\le n\le 7$, then the same conclusion remains
valid for $0<\lambda\le\lambda^{\ast}(1+\chi|\Omega|)^2$ and $\mu_0
=\lambda^{\ast}$.
\end{thm}
\begin{proof}
Let $\{t_i\}_{i=1}^{\infty}$ be a sequence such that $t_i\ge t_0$ for any
$i\ge 1$ and $t_i\to\infty$ as $i\to\infty$. By Corollary 2.4 the sequence
$\{t_i\}_{i=1}^{\infty}$ has a subsequence $\{t_i'\}_{i=1}^{\infty}$ such that
$u(\cdot,t_i')$ converges uniformly in $C^2(\2{\Omega})$ to a solution $v$ 
of ($S_{\lambda}^N$) as $i\to\infty$. Let
\begin{equation}
\lambda_0=\frac{\lambda}{(1+\chi\int_{\Omega}\frac{dy}{1-v(y)})^2}.
\end{equation}
Then $\lambda_0\le\mu\le\mu_0$. Let $\mu_1=\min\{\mu'\ge\lambda_0:w_{\mu'}(x)
\ge v(x)\quad\forall x\in\Omega\}$. Then $\lambda_0\le\mu_1\le\mu_0$ and 
$w_{\mu_1}(x)\ge v(x)$ in $\Omega$.
Let $q(x)=w_{\mu_1}(x)-v(x)$. Then $q(x)\ge 0$ in $\Omega$ and $q(x)=0$ on
$\1\Omega$. Suppose $q(x)\ne 0$ in $\Omega$. Then $\mu_1>\lambda_0$. Since
\begin{equation*}
-\Delta q=\frac{\mu_1}{(1-w_{\mu_1})^2}-\frac{\lambda_0}{(1-v)^2}
>\frac{\lambda_0}{(1-w_{\mu_1})^2}-\frac{\lambda_0}{(1-v)^2}
=\frac{2\lambda_0q}{(1-\xi)^3}\ge 0\quad\mbox{ in }\Omega
\end{equation*}
for some function $\xi(x)$ between $w_{\mu_1}(x)$ and $v(x)$, by the strong 
maximum principle 
\begin{align}
&q(x)>0\quad\mbox{ in }\Omega\quad\mbox{ and }\quad\frac{\1 q}{\1\nu}<0\quad
\mbox{ on }\1\Omega\nonumber\\
\Rightarrow\quad&w_{\mu_1}(x)>v(x)\quad\mbox{ in }\Omega\quad\mbox{ and }
\quad\frac{\1 w_{\mu_1}}{\1\nu}<\frac{\1 v}{\1\nu}\quad\mbox{ on }\1\Omega
\end{align}
where $\1/\1\nu$ is the derivative with respect to the unit exterior normal 
$\nu$ on $\1\Omega$. Let
\begin{equation*}
\3_1=\frac{1}{3}\min_{\1\Omega}\biggl(\frac{\1 v}{\1\nu}
-\frac{\1 w_{\mu_1}}{\1\nu}\biggr).
\end{equation*}
Then $\3_1>0$. Since $\Omega\subset\R^n$ is a smooth bounded domain, 
there exists  $\delta_1>0$ such that for each $x\in\Omega_{\delta_1}$ 
there exists a unique minimizing normalized geodesic $\gamma
=\gamma_x:[0,\rho_1]\to\2{\Omega}$ such that $\gamma(0)=x$, 
$\gamma (\rho_1)\in\1{\Omega}$, $\gamma ([0,\rho_1))
\subset\Omega$, $\gamma'(\rho_1)$ is perpendicular to $\1\Omega$ at
$\gamma (\rho_1)$ where $\rho_1=\mbox{dist}(x,\1\Omega)$ (cf. \cite{Wa}). 
We may also assume that $\delta_1$ is small such that 
\begin{equation}
\frac{\1 }{\1\gamma}(v-w_{\mu_1})\ge\3_1\quad\forall x\in\Omega_{\delta_1}
\end{equation}
where $\1/\1\gamma$ is the partial derivative along the geodesic $\gamma$. 
By Theorem 1.1 and (3.3) there exists $\mu_2\in (\lambda_0,\mu_1)$ such that
\begin{equation}
w_{\mu_2}(x)>v(x)\quad\mbox{ in }\Omega\setminus\Omega_{\delta_1}
\end{equation}
and
\begin{equation}
\left|\frac{\1}{\1\gamma}(w_{\mu_2}-w_{\mu_1})(x)\right|\le\frac{\3_1}{2}
\quad\forall x\in\Omega_{\delta_1}.
\end{equation}
By (3.4) and (3.6), $\forall x\in\Omega_{\delta_1}$,
\begin{equation}
w_{\mu_2}(x)-v(x)=\int_{\rho_1}^0\frac{\1}{\1\gamma}(w_{\mu_2}-v)ds
=\int_0^{\rho_1}\biggl (\frac{\1}{\1\gamma}(v-w_{\mu_1})
+\frac{\1}{\1\gamma}(w_{\mu_1}-w_{\mu_2})\biggr)ds>0.
\end{equation}
By (3.5) and (3.7),
$$
w_{\mu_2}(x)>v(x)\quad\mbox{ in }\Omega.
$$
This contradicts the choice of $\mu_1$. Hence $q(x)\equiv 0$ on $\2{\Omega}$.
Thus 
\begin{equation}
\lambda_0=\mu_1\quad\mbox{ and }\quad v(x)=w_{\mu_1}(x)\quad\mbox{ in }\Omega.
\end{equation}
By (3.2) and (3.8) $\mu_1$ satisfies (1.2). Since $\mu_1$ is uniquely 
determined by (1.2) independent of the subsequence $\{t_i'\}_{i=1}^{\infty}$,
$\{u(\cdot, t_i)\}_{i=1}^{\infty}$ converges uniformly in 
$C^2(\2{\Omega})$ to $w_{\mu_1}$ as $i\to\infty$.
Since the sequence $\{t_i\}_{i=1}^{\infty}$ is arbitrary, $u(\cdot,t)$ 
converges uniformly in $C^2(\2{\Omega})$ to $w_{\mu_1}$ as $i\to\infty$.  

If $1\le n\le 7$, then by Theorem 1.1 and a similar argument as before 
the same conclusion holds for $0<\lambda\le\lambda^{\ast}(1+\chi|\Omega|)^2$ 
and $\mu_0=\lambda^{\ast}$ and the theorem follows.
\end{proof}

\begin{thm}
Let $1\le n\le 7$ and let $\Omega\subset\R^n$ be a smooth convex bounded 
domain. Let $\chi>0$, $0<a_1\le (1+\chi|\Omega|)^2$ and
$$
\3_0=\frac{a_1(1-\|w_{\ast}\|_{\infty})^2}
{(1+\chi\int_{\Omega}\frac{dx}{1-w_{\ast}})^2}.
$$
Then for any $0<\3_1\le\3_0$, $a_1\lambda^{\ast}\le\lambda
\le\lambda^{\ast}(1+\chi\int_{\Omega}\frac{dx}{1-\3_1w_{\ast}})^2$, 
and $u_0\in L^1(\Omega)$ such that $2\3_1w_{\ast}\le u_0\le w_{\ast}$
in $\Omega$, there exists a unique global solution of ($P_{\lambda}$) 
satisfying
\begin{equation}
\3_1 w_{\ast}\le u\le w_{\ast}\quad\mbox{ in }\Omega\times (0,\infty).
\end{equation}
Moreover $u(\cdot,t)$ converges uniformly in $C^2(\2{\Omega})$ to 
$w_{\mu_1}$ as $t\to\infty$ where $\mu_1>0$ is uniquely given by (1.2).
\end{thm}
\begin{proof}
Note that uniqueness of solution of ($P_{\lambda}$) follows by Theorem 2.1. 
We next prove the existence of global solution of 
($P_{\lambda}$). We divide the proof into two cases.

\noindent $\underline{\text{\bf Case 1}}$: $u_0\in C^{\infty}(\2{\Omega})$.

By Theorem 2.2 there exists $T'>0$ such that ($P_{\lambda}$) has a 
non-negative solution $u$ in $\Omega\times (0,T')$. By the parabolic \
Schauder estimates \cite{LSU} $u\in C^{2+\beta,1+(\beta/2)}
(\2{\Omega}\times [0,T'))$ for some constant $0<\beta<1$.  
Since $\Delta w_{\ast}<0$ in $\Omega$, $w_{\ast}>0$ in $\Omega$ and 
$w_{\ast}=0$ on $\1\Omega$, by the Hopf Lemma,
\begin{equation}
\frac{\1 w_{\ast}}{\1\nu}<0\quad\mbox{ on }\1\Omega\quad
\end{equation}
where $\nu$ the unit outward normal on $\1\Omega$. Since $u_0-2\3_1
w_{\ast}\ge 0$ in $\Omega$ and $u_0-2\3_1 w_{\ast}=0$ on $\1\Omega$, 
by (3.10), 
\begin{equation}
\frac{\1 u_0}{\1\nu}\le 2\3_1\frac{\1 w_{\ast}}{\1\nu}
<\3_1\frac{\1 w_{\ast}}{\1\nu}\quad\mbox{ on }\1\Omega.
\end{equation}
Let
$$
\3_2=\frac{1}{3}\min_{\1\Omega}\left(\3_1\frac{\1 w_{\ast}}{\1\nu}
-\frac{\1 u_0}{\1\nu}\right).
$$ 
By (3.11), $\3_2>0$. Then similar to the proof of Theorem 3.5 there exists
$\delta_1>0$ such that for each $x\in\Omega_{\delta_1}$ there 
exists a unique minimizing normalized geodesic $\gamma=\gamma_x:
[0,\rho_1]\to\2{\Omega}$ 
such that $\gamma(0)=x$, $\gamma (\rho_1)\in\1{\Omega}$, $\gamma ([0,\rho_1))
\subset\Omega$, $\gamma'(\rho_1)$ is perpendicular to $\1\Omega$ at
$\gamma (\rho_1)$ where $\rho_1=\mbox{dist}(x,\1\Omega)$. We may also assume
that $\delta_1$ is small such that 
\begin{equation}
\frac{\1 }{\1\gamma}(\3_1w_{\ast}-u_0)\ge\3_2\quad\forall x\in\Omega_{\delta_1}
\end{equation}
where $\1/\1\gamma$ is the partial derivative along the geodesic $\gamma$.
Let
\begin{equation}
\3_3=\min_{\Omega\setminus\Omega_{\delta_1}}(u_0-\3_1w_{\ast}).
\end{equation}
Then $\3_3>0$. Since $u\in C^{2+\beta,1+(\beta/2)}(\2{\Omega}
\times [0,T'))$, there exists $0<T_1<T'$ such that
\begin{equation}
\|u(\cdot,t)-u_0\|_{L^{\infty}(\Omega)}<\frac{\3_3}{2}\quad\forall
0\le t\le T_1
\end{equation}
and
\begin{equation}
\left|\frac{\1 }{\1\gamma}(u(x,t)-u_0(x))\right|\le\frac{\3_2}{2}\quad\forall 
x\in\Omega_{\delta_1}, 0\le t\le T_1.
\end{equation}
By (3.12), (3.13), (3.14) and (3.15),
\begin{equation}
u\ge\3_1w_{\ast}\quad\mbox{ in }(\Omega\setminus\Omega_{\delta_1})
\times [0,T_1]
\end{equation}
and
\begin{equation}
\frac{\1 }{\1\gamma}(\3_1w_{\ast}(x)-u(x,t))\ge 0\quad\forall x\in
\Omega_{\delta_1},0\le t\le T_1.
\end{equation}
Let $x\in\Omega_{\delta_1}$ and $0\le t\le T_1$. Then by (3.17),
\begin{equation}
\3_1w_{\ast}(x)-u(x,t)=\int_{\rho_1}^0\frac{\1 }{\1\gamma}(\3_1w_{\ast}
(\gamma (s))-u(\gamma (s),t))\,ds\le 0\quad\forall x\in\Omega_{\delta_1},
0\le t\le T_1.
\end{equation}
By (3.16) and (3.18),
\begin{equation}
u(x,t)\ge\3_1 w_{\ast}(x)\quad\forall x\in\Omega,0\le t\le T_1.
\end{equation}
Let $T=\sup\{T_2>0:u\ge\3_1w_{\ast}\mbox{ in }\Omega\times [0,T_2]\}$. 
Then $T\ge T_1$ and
\begin{align}
u_t\le&\Delta u 
+\frac{\lambda}{(1-u)^2(1+\chi\int_{\Omega}\frac{dy}{1-\3_1w_{\ast}})^2}
\quad\mbox{ in }\Omega\times (0,T)\nonumber\\
\le&\Delta u+\frac{\lambda^\ast}{(1-u)^2}\qquad\qquad\qquad\qquad
\quad\mbox{in }\Omega\times (0,T).
\end{align}
By (3.20) and Lemma 3.3,
\begin{equation}
u\le w_{\ast}\quad\mbox{in }\Omega\times (0,T).
\end{equation}
By the Schauder estimates we can extend $u$ to a function in
$C^{2+\beta,1+(\beta/2)}(\2{\Omega}\times [0,T])$.
By (3.21) and ($P_{\lambda}$),
\begin{equation}
u_t\ge\Delta u
+\frac{\lambda}{(1-u)^2(1+\chi\int_{\Omega}\frac{dy}{1-w_{\ast}})^2}
\quad\mbox{ in }\Omega\times (0,T).
\end{equation}
Let $\psi=\3_1w_{\ast}$. Then
\begin{align}
\Delta\psi+\frac{\lambda}{(1-\psi)^2(1+\chi\int_{\Omega}
\frac{dx}{1-w_{\ast}})^2}=&-\frac{\3_1\lambda^{\ast}}{(1-w_{\ast})^2}
+\frac{\lambda}{(1-\psi)^2(1+\chi\int_{\Omega}\frac{dx}{1-w_{\ast}})^2}
\nonumber\\
\ge&\lambda^{\ast}\left(\frac{a_1}{(1+\chi\int_{\Omega}
\frac{dx}{1-w_{\ast}})^2}-\frac{\3_1}{(1-w_{\ast})^2}\right)\nonumber\\
\ge&0\quad\mbox{ in }\Omega\times (0,T).
\end{align}
By (3.19), (3.22) and (3.23),
\begin{equation}
\left\{\begin{aligned}
&(u-\psi)_t\ge\Delta (u-\psi)+\lambda F(u,\psi)(u-\psi)
\ge\Delta (u-\psi)+\lambda a_0(u-\psi)
\quad\mbox{ in }\Omega\times (0,T)\\
&u(x,0)-\psi(x)>\3_1w_{\ast}>0\quad\mbox{ in }\Omega
\end{aligned}\right.
\end{equation}
where
$$
F(u,\psi)=\frac{2-u-\psi}{(1-u)^2(1-\psi)^2(1+\chi\int_{\Omega}
\frac{dx}{1-w_{\ast}})^2}
$$
and
$$
a_0=\min_{\2{\Omega}\times [0,T]}F(u,\psi)>0.
$$
Hence
\begin{equation}
\left\{\begin{aligned}
&(e^{-\lambda a_0t}(u-\psi))_t\ge\Delta (e^{-\lambda a_0t}(u-\psi))
\quad\mbox{ in }\Omega\times (0,T]\\
&u(x,0)-\psi(x)>0\quad\mbox{ in }\Omega
\end{aligned}\right.
\end{equation}
By (3.25) and the strong maximum principle,
\begin{equation}
u-\psi>0\quad\mbox{ in }\Omega\times (0,T]\quad\Rightarrow\quad
u(x,T)>\3_1 w_{\ast}(x)\quad\mbox{ in }\Omega
\end{equation}
and
\begin{equation}
\frac{\1}{\1\nu}(u(x,T)-\3_1w_{\ast}(x))<0\quad\mbox{ on }\1\Omega
\end{equation}
By (3.26) and (3.27) and an argument similar to the one before there exists
a constant $\delta>0$ such that
\begin{equation}
u(x,T)\ge(\3_1+\delta)w_{\ast}(x)\quad\mbox{ in }\Omega
\end{equation}
By repeating the above argument there exists a constant $T_2>0$ such that
there exists a solution $\4{u}$ of ($P_{\lambda}$) in $\Omega\times (0,T_2)$
with initial value $u(x,T)$ such that $\4{u}\ge\3_1 w_{\ast}$ in 
$\Omega\times (0,T_2)$. Let $u(x,t)=\4{u}(x,t-T)$ for all $x\in\Omega$,
$T\le t\le T+T_2$. Then $u$ is a solution of ($P_{\lambda}$) in 
$\Omega\times (0,T+T_2)$ such that $u\ge\3_1 w_{\ast}$ in 
$\Omega\times (0,T_2)$. This contradicts the maximality of $T$. Hence 
$T=\infty$.

\noindent $\underline{\text{\bf Case 2}}$: $u_0\in L^1(\Omega)$

We choose a sequence of function $u_{0,k}\in C^{\infty}(\Omega)$ satisfying
$$
2\3_1 w_{\ast}\le u_{0,k}\le w_{\ast}\quad\mbox{ in }\Omega
$$
such that $u_{0,k}\to u_0$ in $L^1(\Omega)$ as $k\to\infty$.
For each $k\ge 1$ by case 1 there exists a unique global solution $u_k$ of 
($P_{\lambda}$) with initial value $u_{0,k}$ satisfying 
$$
\3_1 w_{\ast}\le u_k\le w_{\ast}\quad\mbox{ in }\Omega\times (0,\infty).
$$  
By the parabolic Schauder estimates \cite{LSU} the sequence 
$\{u_k\}_{k=1}^{\infty}$ are uniformly bounded in $u\in 
C^{2+\beta,1+(\beta/2)}(\2{\Omega}\times (\delta,\infty))$ for some
constant $0<\beta<1$ and any $\delta>0$. 
By the Ascoli-Arzel\'a theorem and a 
diagonalization argument the sequence $\{u_k\}_{k=1}^{\infty}$ has a 
subsequence which we may assume without loss of generality to be the sequence
itself that converges uniformly in $C^{2+\beta,1+(\beta/2)}(\2{\Omega}
\times (\delta,1/\delta))$ for any $0<\delta<1$ to some function $u$ 
as $k\to\infty$. Then $u$  satisfies
\begin{equation*}\left\{\begin{aligned}
\frac{\1 u}{\1 t}=&\Delta u
+\frac{\lambda}{(1-u)^2(1+\chi\int_\Omega\frac{dy}{1-u(y,t)})^2} 
\quad\quad\mbox{in }\Omega\times (0,\infty)\\
u=&0\qquad\qquad\qquad \qquad\qquad\quad\quad\quad\quad\quad
\mbox{on }\partial\Omega\times (0,\infty)
\end{aligned}\right.
\end{equation*} 
and (3.9). Since each $u_k$ satisfies (2.30) in $\Omega\times (0,\infty)$, 
letting $k\to\infty$ we get that $u$ satisfies (2.14) in $\Omega\times 
(0,\infty)$. Letting $t\to 0$ in (2.14), by (3.9) $u(\cdot,t)\to u_0$ in 
$L^1(\Omega)$ as $t\to\infty$. Hence $u$ is the global solution of 
($P_{\lambda}$). 

By (3.9) and an argument similar to the proof of Theorem 3.5
$u(\cdot,t)$ converges uniformly in $C^2(\2{\Omega})$ to 
$w_{\mu_1}$ as $t\to\infty$ where $\mu_1>0$ is uniquely given by
(1.2) and the theorem follows.
\end{proof}

\begin{thm}
Let $1\le n\le 7$ and let $\Omega\subset\R^n$ be a smooth convex bounded 
domain. Let $0<\delta<1/2$, $\chi>0$ and $0<\lambda_2\le\lambda^{\ast}$ 
satisfy
$$
\lambda_2\frac{(1-2\delta)}{(1-\|w_{\ast}\|_{\infty})^2}
\left(1+\chi\int_{\Omega}\frac{dx}{1-w_{\ast}}\right)^2\le\lambda
\le\lambda^{\ast}\left(1+\chi\int_{\Omega}\frac{dx}{1-(1-2\delta)
w_{\lambda_2}}\right)^2.
$$
Let $u_0\in L^1(\Omega)$ satisfy $(1-\delta)w_{\lambda_2}\le u_0
\le w_{\ast}$ in $\Omega$.
Then there exists a unique global solution of ($P_{\lambda}$) 
satisfying
\begin{equation}
(1-2\delta)w_{\lambda_2}\le u\le w_{\ast}\quad\mbox{ in }\Omega\times 
(0,\infty). 
\end{equation}
Moreover $u(\cdot,t)$ converges uniformly in $C^2(\2{\Omega})$ to 
$w_{\mu_1}$ as $t\to\infty$ where $\mu_1>0$ is uniquely given by (1.2).
\end{thm}
\begin{proof}
By an approximation argument similar to the proof of Theorem 3.6 it suffices
to prove the existence of global solution of ($P_{\lambda}$) for the case 
$u_0\in C^{\infty}(\2{\Omega})$. By Theorem 2.2 and an argument similar 
to the proof of Theorem 3.6 there exists a maximal time $T>0$ such that 
($P_{\lambda}$) has a solution $u$ in $\Omega\times (0,T)$ which satisfies
\begin{equation}
u\ge (1-2\delta)w_{\lambda_2}\quad\mbox{ in }\Omega\times (0,T).
\end{equation}
By ($P_{\lambda}$) and (3.30),
\begin{align}
u_t\le&\Delta u+\frac{\lambda}{(1-u)^2(1+\chi\int_{\Omega}
\frac{dy}{(1-(1-2\delta)w_{\lambda_2})})^2}
\quad\mbox{ in }\Omega\times (0,T)\nonumber\\
\le&\Delta u+\frac{\lambda^\ast}{(1-u)^2}\qquad\qquad\qquad\qquad
\quad\mbox{in }\Omega\times (0,T).
\end{align}
By (3.31) and Lemma 3.3,
\begin{equation}
u\le w_{\ast}\quad\mbox{ in }\Omega\times (0,T).
\end{equation}
By (3.30), (3.32) and the parabolic Schauder estimates \cite{LSU}, 
$u\in C^{2+\beta,1+(\beta/2)}(\2{\Omega}\times [0,T))$ for some constant
$0<\beta<1$. Hence
we can extend $u$ to a function in $C^{2+\beta,1+(\beta/2)}(\2{\Omega}
\times [0,T])$. Then by ($P_{\lambda}$) and (3.32),
\begin{equation}
u_t\ge\Delta u+\frac{\lambda}{(1-u)^2(1+\chi\int_{\Omega}
\frac{dy}{(1-w_{\ast})})^2}
\quad\mbox{ in }\Omega\times (0,T).
\end{equation}
Let $\psi=(1-2\delta)w_{\lambda_2}$. Then
\begin{align}
\Delta\psi+
\frac{\lambda}{(1-\psi)^2(1+\chi\int_{\Omega}\frac{dx}{1-w_{\ast}})^2}
=&-\lambda_2\frac{(1-2\delta)}{(1-w_{\lambda_2})^2}+
\frac{\lambda}{(1-\psi)^2(1+\chi\int_{\Omega}\frac{dx}{1-w_{\ast}})^2}
\nonumber\\
\ge&\frac{1}{(1-w_{\lambda_2})^2}\left(
\frac{\lambda (1-w_{\lambda_2})^2}
{(1+\chi\int_{\Omega}\frac{dx}{1-w_{\ast}})^2}
-\lambda_2(1-2\delta)\right)\nonumber\\
\ge&0.
\end{align}
By (3.33), (3.34) and an argument similar to the proof of Theorem 3.6,
\begin{equation}
u-\psi>0\quad\mbox{ in }\Omega\times (0,T]\quad\Rightarrow\quad
u(x,T)>(1-2\delta)w_{\lambda_2}(x)\quad\mbox{ in }\Omega
\end{equation}
and
\begin{equation}
\frac{\1}{\1\nu}(u(x,T)-(1-2\delta)w_{\lambda_2}(x))<0\quad\mbox{ on }
\1\Omega.
\end{equation}
By (3.35) and (3.36) and an argument similar to the proof of Theorem 3.6,
there exists a constant $\delta_1>0$ such that
\begin{equation}
u(x,T)\ge (1-2\delta+\delta_1)w_{\lambda_2}(x)\quad\mbox{ on }\Omega.
\end{equation}
Then similar to the proof of Theorem 3.6 by (3.37) $u$ can be extended to 
a solution of ($P_{\lambda}$) in $\Omega\times (0,T+T_1)$ for some 
$T_1>0$ such that
\begin{equation*}
u\ge (1-2\delta)w_{\lambda_2}\quad\mbox{ in }\Omega\times (0,T+T_1).
\end{equation*}
This contradicts the maximality of $T$. Hence $T=\infty$. By an argument 
similar to the proof of Theorem 3.6 $u(\cdot,t)$ converges uniformly in 
$C^2(\2{\Omega})$ to $w_{\mu_1}$ as $t\to\infty$ where $\mu_1>0$ 
is uniquely given by (1.2) and the theorem follows.
\end{proof}

By an argument similar to the proof of Theorem 3.6 and Theorem 3.7 
we have the following result.

\begin{thm}
Let $1\le n\le 7$ and $\Omega\subset\R^n$ be a smooth convex bounded 
domain. Let $0<\delta<1/2$, $\chi>0$, $0<\lambda\le\lambda^{\ast}
(1+\chi|\Omega|)^2$, 
$$
\mu=\frac{\lambda}{(1+\chi|\Omega|)^2}\quad\mbox{ and }
\quad\mu'=\frac{\lambda}{(1+\chi\int_{\Omega}
\frac{dx}{1-w_{\ast}})^2}.
$$
Let $u_0\in L^1(\Omega)$ satisfy $w_{(1-\delta)\mu'}
\le u_0\le w_{\mu}$ in $\Omega$. Then there exists a unique global 
solution of ($P_{\lambda}$) satisfying
\begin{equation*}
w_{(1-2\delta)\mu'}\le u\le w_{\mu}\quad\mbox{ in }\Omega\times (0,\infty).
\end{equation*}
\end{thm}

\section{Quenching behaviour}
\setcounter{equation}{0}

In this section we will prove the quenching behaviour of the solution 
of ($P_{\lambda}$) when $\lambda$ is large. 
We first start with a technical lemma.

\begin{lem}
Let $u_0\equiv 0$ on $\Omega$ and let $\chi$ satisfy
\begin{equation}
0<\chi<\frac{1}{|\Omega|}.
\end{equation}
Then there exist constants $\lambda_1>0$ and $C_{\lambda_1}$ such that for 
any $\lambda\ge\lambda_1$ and any global solution $u$ of 
$(P_{\lambda})$ there exists a  sequence $\{t_i\}_{i=1}^{\infty}$, 
$t_i\to\infty$ as $i\to\infty$, such that
\begin{equation}
\int_{\Omega}u_t^2(x,t_i)\,dx\to 0\quad\mbox{ as }i\to\infty  
\end{equation}
and
\begin{equation}
\int_{\Omega}\frac{dx}{(1-u(x,t_i))^2}\le C_{\lambda_1}
\quad\forall i\in \mathbb{Z}^+.
\end{equation}
\end{lem}
\begin{proof}
Let $0<\delta<1$. By (4.1) we can choose constants $0<\3<1$, $\lambda_1>0$ 
and $K>0$ such that
\begin{equation}
a_0=\frac{1+\3}{1-\3}\biggl(
\frac{\delta\chi^2\sqrt{|\Omega|}}{\lambda_1}+
\frac{2\chi|\Omega|}{1+\chi|\Omega|}\biggr )(1+K^{-1})<1.
\end{equation}
Let $\lambda\ge\lambda_1$ and let $u$ be a global solution 
of $(P_{\lambda})$. By Theorem 2.1 and Theorem 2.2, $u\ge 0$ on 
$\Omega\times (0,\infty)$. By (2.23) of Theorem 2.3 there exists a sequence 
$\{t_i\}_{i=1}^{\infty}$, $t_i\to\infty$ as $i\to\infty$, such that 
(4.2) holds. Let $u_i=u(x,t_i)$. By (4.2) we can assume without loss
of generality that 
\begin{equation}
\int_{\Omega}u_{i,t}^2\,dx\le\frac{\delta^2}{|\Omega|}\quad\forall
i\in\Z^+.
\end{equation}
Multiplying ($P_{\lambda}$) by $u$, integrating over 
$\Omega$ and putting $t=t_i$,
\begin{equation}
\int_{\Omega}u_iu_{i,t}\,dx+\int_{\Omega}|\nabla u_i|^2\,dx=\lambda
\frac{\int_{\Omega}\frac{u_i}{(1-u_i)^2}\,dx}{(1+\chi\int_{\Omega}
\frac{dx}{1-u_i})^2}.
\end{equation}
Now
\begin{align}
\int_{\Omega}\frac{u_i}{(1-u_i)^2}\,dx=&\int_{\Omega}\frac{dx}{(1-u_i)^2}
-\int_{\Omega}\frac{dx}{1-u_i}\nonumber\\
\ge&\int_{\Omega}\frac{dx}{(1-u_i)^2}-|\Omega|^{\frac{1}{2}}\biggl (
\int_{\Omega}\frac{dx}{(1-u_i)^2}\biggr )^{\frac{1}{2}}\nonumber\\
\ge&(1-\3)\int_{\Omega}\frac{dx}{(1-u_i)^2}-\frac{|\Omega|}{4\3}
\end{align}
and
\begin{align}
\biggl (1+\chi\int_{\Omega}\frac{dx}{1-u_i}\biggr )^2
=&1+2\chi\int_{\Omega}\frac{dx}{1-u_i}+\chi^2\biggl (
\int_{\Omega}\frac{dx}{1-u_i}\biggr )^2\nonumber\\
\le&1+2\chi|\Omega|^{\frac{1}{2}}\biggl (
\int_{\Omega}\frac{dx}{(1-u_i)^2}\biggr )^{\frac{1}{2}}
+\chi^2|\Omega|\int_{\Omega}\frac{dx}{(1-u_i)^2}\nonumber\\
\le&1+\frac{1}{\3}+(1+\3)\chi^2|\Omega|\int_{\Omega}\frac{dx}{(1-u_i)^2}.
\end{align}
If there exists a subsequence of $u_i$ which we may assume without loss of 
generality to be the sequence itself such that 
$$
(1+\3)\chi^2|\Omega|\int_{\Omega}\frac{dx}{(1-u_i)^2}\le (1+\3^{-1})K
\quad\forall i\in\mathbb{Z}^+,
$$
then (4.3) follows and we are done. Suppose no such subsequence exists. 
Then there exists a subsequence of $u_i$ which we may assume without loss of 
generality to be the sequence itself such that 
\begin{equation}
(1+\3)\chi^2|\Omega|\int_{\Omega}\frac{dx}{(1-u_i)^2}>(1+\3^{-1})K
\quad\forall i\in\mathbb{Z}^+.
\end{equation}
By (4.1), (4.2), (4.5), (4.6), (4.7), (4.8), (4.9) and Theorem 2.3,
\begin{align}
&\lambda\biggl\{(1-\3)\int_{\Omega}\frac{dx}{(1-u_i)^2}-\frac{|\Omega|}{4\3}
\biggr\}\nonumber\\
\le&\biggl(\int_{\Omega}u_iu_{i,t}\,dx+\int_{\Omega}|\nabla u_i|^2\,dx
\biggr)\biggl(1+\3^{-1}+(1+\3)\chi^2|\Omega|\int_{\Omega}\frac{dx}{(1-u_i)^2}
\biggr)\nonumber\\
\le&(1+K^{-1})(1+\3)\chi^2|\Omega|
\biggl(\biggl (\int_{\Omega}u_i^2\biggr )^\frac{1}{2}\biggl (
\int_{\Omega}u_{i,t}^2\biggr )^\frac{1}{2}
+\frac{2\lambda}{\chi(1+\chi|\Omega|)}\biggr)
\int_{\Omega}\frac{dx}{(1-u_i)^2}\nonumber\\
\le&(1+K^{-1})(1+\3)\lambda
\biggl(\frac{\delta\chi^2\sqrt{|\Omega|}}{\lambda_1}
+\frac{2\chi|\Omega|}{1+\chi|\Omega|}\biggr)
\int_{\Omega}\frac{dx}{(1-u_i)^2}
\end{align}
Hence by (4.4) and (4.10),
$$
\int_{\Omega}\frac{dx}{(1-u_i)^2}\le\frac{|\Omega|}{4\3(1-\3)(1-a_0)}\quad
\forall i\in\mathbb{Z}^+.
$$
and the lemma follows.
\end{proof}

\begin{thm}
Let $u_0\equiv 0$ on $\Omega$ and let $\chi$ satisfy (4.1). Let 
$\lambda_1>0$ be given by Lemma 4.1. Then there exists a constant 
$\lambda_2\ge\lambda_1$ such that for any $\lambda>\lambda_2$ and any \
solution $u$ of ($P_{\lambda}$) there exists $T_1>0$ such that
\begin{equation}
\lim_{t\nearrow T_1}\sup_{\tiny\begin{array}{c}
x\in\Omega\end{array}}u(x,t)=1.
\end{equation}
\end{thm}
\begin{proof}
Let $\lambda_1$ and $C_{\lambda_1}$ be as in Lemma 4.1.
Let $\lambda>\lambda_2$ for some constant 
$\lambda_2\ge\lambda_1$ to be determined later. Suppose $u$ is a global 
solution of ($P_{\lambda}$). By Theorem 2.1 and Theorem 2.2, $u\ge 0$
in $\Omega\times (0,\infty)$. By Lemma 4.1 there exists a sequence 
$\{t_i\}_{i=1}^{\infty}$, $t_i\to\infty$ as $i\to\infty$, such that
(4.2) and (4.3) holds.
Let $u_i=u(x,t_i)$. Since $0\le u_i<1$, by (4.3) and Theorem 2.3 
there exist  $u_{\infty}\in H_0^1(\Omega)$, $0\le u_{\infty}\le 1$
in $\Omega$, $0\le g\in L^2(\Omega)$ such that  
\begin{equation*}
\int_{\Omega}g^2\,dx\le C_{\lambda_1}
\end{equation*}
and a subsequence of $\{u_i\}$ 
which we may assume without loss of generality to be the sequence 
$\{u_i\}$ itself such that
\begin{equation}
\left\{\begin{aligned}
&u_i\to u_{\infty}\qquad\qquad\mbox{ weakly in }H_0^1(\Omega)
\quad\mbox{ as }i\to\infty\\
&u_i\to u_{\infty}\qquad\qquad\mbox{ weakly in }L^2(\Omega)
\quad\mbox{ as }i\to\infty\\
&(1-u_i)^{-1}\to g\quad\mbox{ weakly in }L^2(\Omega)
\quad\mbox{ as }i\to\infty.
\end{aligned}\right. 
\end{equation}
Letting $i\to\infty$ in (4.3), by the Fatou Lemma,
\begin{equation}
\int_{\Omega}\frac{dx}{(1-u_{\infty})^2}\le C_{\lambda_1}.
\end{equation}
Hence $u_{\infty}(x)<1$ a.e. $x\in\Omega$.
Let $0\le\eta\in C_0^{\infty}(\Omega)$. Multiplying ($P_{\lambda}$) 
by $\eta$, integrating over $\Omega$ and putting $t=t_i$,
\begin{equation}
\int_{\Omega}u_{i,t}\eta\,dx=-\int_{\Omega}\nabla u_i\cdot\nabla\eta\,dx
+\lambda\frac{\int_{\Omega}\frac{\eta}{(1-u_i)^2}\,dx}{(1+\chi\int_{\Omega}
\frac{dx}{1-u_i})^2}\quad\forall i\in\mathbb{Z}^+. 
\end{equation}
Letting $i\to\infty$ in (4.14), by (4.2) and (4.12),
\begin{equation}
-\int_{\Omega}\nabla u_{\infty}\cdot\nabla\eta\,dx
+\lambda\frac{\int_{\Omega}\frac{\eta}{(1-u_{\infty})^2}\,dx}
{(1+\chi\int_{\Omega}g\,dx)^2}\le C\lim_{i\to\infty}
\|u_{i,t}\|_{L^2(\Omega)}=0.
\end{equation}
Putting $\3=1$ in (4.8) and letting $i\to\infty$, by (4.3) and (4.12)
we have
\begin{equation}
(1+\chi\int_{\Omega}g\,dx)^2
\le 2+2\chi^2|\Omega|C_{\lambda_1}=K_1\quad (\mbox{say}).
\end{equation}
Let $\lambda_2=\max (\lambda_1,\lambda^{\ast}K_1)$. Then by (4.15)
and (4.16),
\begin{equation*}
-\int_{\Omega}\nabla u_{\infty}\cdot\nabla\eta\,dx
+\frac{\lambda}{K_1}\int_{\Omega}\frac{\eta}{(1-u_{\infty})^2}\,dx\le 0.
\end{equation*}
Hence $u_{\infty}$ is a weak supersolution of ($S_{\lambda/K_1}$) and
$\lambda/K_1>\lambda^{\ast}$. Let $\lambda^{\ast}<\lambda_3<\lambda/K_1$. 
By an argument similar to the proof of Proposition 5.3 of \cite{GG1} there 
exists a classical solution of ($S_{\lambda_3}$). This contradicts the 
maximality of $\lambda^{\ast}$. Hence there exists $T_1>0$ such that 
(4.11) holds and the theorem follows.
\end{proof}

\begin{thm}
Let $\Omega=B_R$ and $0\le u_0\in L^1(B_R)$ be a radially symmetric
monotone decreasing function which satisfies (0.1) for some constant 
$0<a<1$. Let $\chi>0$. Then there exists a constant $C_3>0$ and such 
that for any $\lambda>\lambda_0=C_3\mu_1$ and any solution $u$ of 
($P_{\lambda})$, $u$ quenches in a finite time 
\begin{equation}
T_{\lambda}\le\frac{C_3}{\lambda-\lambda_0}.
\end{equation}
\end{thm}
\begin{proof}
Let $u$ be a global solution of ($P_{\lambda})$ and let
$$
E(t)=\int_{\Omega}u(x,t)\phi_1(x)\,dx.
$$ 
By Theorem 2.1 and 
Theorem 2.2 $u\ge 0$ and $u$ is radially symmetric in $\Omega\times 
(0,\infty)$. Hence $u(r,t)=u(|x|,t)$ where $r=|x|$. Since $u_0(r)$ is a 
monotone decreasing function of $0<r<R$, by the strong maximum principle 
and an argument similar to the proof of Theorem 1.5 of \cite{Hs} 
$u_r(r,t)<0$ for all $0<r<R$ and $t>0$. Then by an argument similar to 
the proof of Proposition 1.3 there exist constants $C_1>0$ and $C_2>0$ 
such that (1.4) and (1.6) hold with $v$ being replaced by $u(x,t)$.

Multiplying ($P_{\lambda}$) by $\phi_1$ and integrating over $\Omega$, by the 
Green theorem, (1.4), and (1.6),
\begin{align}
\frac{d}{dt}E(t)=&\frac{d}{dt}\left(\int_{\Omega}u\phi_1\,dx\right)
=\int_{\Omega}\phi_1\Delta u\,dx
+\lambda\frac{\int_{\Omega}\frac{\phi_1}{(1-u)^2}\,dx}
{(1+\chi\int_{\Omega}\frac{dx}{1-u})^2}\nonumber\\
\ge&-\mu_1\int_{\Omega}u\phi_1\,dx
+\lambda\frac{\int_{\Omega}\frac{\phi_1}{(1-u)^2}\,dx}
{2+C_2\chi^2\int_{\Omega}\frac{\phi_1}{(1-u)^2}\,dx}
\nonumber\\
\ge&-\mu_1 E(t)+\lambda\frac{\int_{\Omega}\phi_1\,dx}
{2+C_2\chi^2\int_{\Omega}\phi_1\,dx}\nonumber\\
\ge&-\mu_1+\frac{\lambda}{C_3}\nonumber\\
\ge&\frac{\lambda-\lambda_0}{C_3}
\end{align}
where $C_3=2+C_2\chi^2$ and $\lambda_0=C_3\mu_1$. By (4.18),
\begin{equation*}
\frac{(\lambda-\lambda_0)}{C_3}t\le E(t)\le 1.
\end{equation*}
Hence $u$ quenches in a finite time $T_{\lambda}$ which satisfies (4.17) and
the theorem follows.
\end{proof}

$$
\mbox{Acknowledgements}
$$

I would like to thank the referee for the detail comments and suggestions
on the paper which leads to great improvement for the paper.

\end{document}